\documentclass[twoside,draft]{article}

\usepackage{amsfonts}
\usepackage{stmaryrd}
\usepackage{amssymb}
\usepackage{euscript}
\usepackage{amsthm}
\usepackage{amsmath}
\usepackage{amscd}
\usepackage{latexsym}
\usepackage{mathrsfs}
\usepackage{graphicx}
\usepackage{color}
\usepackage{dsfont}
\usepackage{bm}

\usepackage{mathtools}
\usepackage{nccmath}
\usepackage{amsfonts}
\usepackage{stmaryrd}
\usepackage{amssymb}
\usepackage{euscript}
\usepackage{amsthm}
\usepackage{amsmath}
\usepackage{amscd}
\usepackage{latexsym}
\usepackage{mathrsfs}
\usepackage{graphicx}
\usepackage{pgf}
\usepackage{color}
\usepackage{xcolor}
\usepackage{dsfont}
\usepackage{bm}
\usepackage{cleveref}
\usepackage{enumitem}%
\usepackage{tikz}
\usetikzlibrary{decorations.pathreplacing}

\newtheorem{theorem}{Theorem}[section]
\newtheorem{lem}{Lemma}[section]
\newtheorem{pro}{Proposition}[section]
\newtheorem{cor}{Corollary}[section]
\newtheorem{rem}{Remark}[section]
\newtheorem{rems}{Remarks}[section]
\newtheorem{ex}{Example}[section]
\newtheorem{defi}{Definition}[section]
\newtheorem{hyp}{Assumption}[section]

\newcommand{\bt}{\begin{theorem}}
\newcommand{\et}{\end{theorem}}
\newcommand{\bl}{\begin{lem}}
\newcommand{\el}{\end{lem}}
\newcommand{\bp}{\begin{pro}}
\newcommand{\ep}{\end{pro}}
\newcommand{\bcor}{\begin{cor}}
\newcommand{\ecor}{\end{cor}}
\newcommand{\lab }{\label }
\newcommand{\bd}{\begin{defi} \rm }
\newcommand{\ed}{\end{defi}}
\newcommand{\brem }{\begin{rem} \rm }
\newcommand{\erem }{\end{rem}}
\newcommand{\brems }{\begin{rems} \rm }
\newcommand{\erems }{\end{rems}}
\newcommand{\bhyp }{\begin{hyp} \rm }
\newcommand{\ehyp }{\end{hyp}}
\newcommand{\bex}{\begin{ex} \rm }
\newcommand{\eex}{\end{ex}}
\numberwithin{equation}{section}
\theoremstyle{plain}

\newcommand{\Nogg}{N^{o,\GG}}

\def\wt{\widetilde}
\def\wh{\widehat }
\def\I{\mathds{1}}

\newcommand{\DB}{D}

\newcommand{\cOff}{\mathcal{O}(\FF)}
\newcommand{\cOgg}{\mathcal{O}(\GG)}
\newcommand{\cPff}{\mathcal{P}(\FF)}
\newcommand{\cPgg}{\mathcal{P}(\GG)}
\newcommand{\cOffd}{\mathcal{O}_d(\FF)}

\newcommand{\cPffd}{\mathcal{P}_d(\FF)}
\newcommand{\cPggd}{\mathcal{P}_d(\GG)}
\newcommand{\cPrff}{\mathcal{P}r(\FF)}

\newcommand{\cPrffd}{\mathcal{P}r\!_d(\FF)}

\newcommand{\Martg}{\mathcal{M}(\GG)}

\newcommand{\Martf}{\mathcal{M}(\FF)}

\newcommand{\Martft}{\mathcal{M}^{\xtheta}(\FF)}

\newcommand{\Martfl}{\mathcal{M}_{loc}(\FF)}
\newcommand{\Martglt}{\mathcal{M}^{\xtheta}_{loc}(\GG)}
\newcommand{\Martflt}{\mathcal{M}^{\xtheta}_{loc}(\FF)}

\newcommand{\whJ}{\wh{J}}

\newcommand{\llb}{\llbracket}
\newcommand{\rrb}{\rrbracket}

\newcommand{\seq}[1]{{\lbrace #1 \rbrace}}

\newcommand{\cEgg}{\wh{\mathcal{E}}}

\newcommand{\xtheta }{\vartheta }
\newcommand{\tgg}{\wh{\tau}}
\newcommand{\sgg}{\wh{\sigma}}

\newcommand{\xN}{K}
\newcommand{\xNX}{K}

\newcommand{\Ygg}{\wh \xY}
\newcommand{\Xgg}{\wh X}
\newcommand{\Xff}{X}

\newcommand{\Zgg}{\wh \xZ}
\newcommand{\Ugg}{\wh \xU}

\newcommand{\Lgg}{\wh \xL}

\newcommand{\ygg}{\wh \xy}
\newcommand{\ugg}{\wh \xu}
\newcommand{\zgg}{\wh \xz}
\newcommand{\lgg}{\wh \xl}

\newcommand{\cT}{{\cal T}}

\newcommand{\cTg}{\mathcal{\wh T}}

\newcommand{\wtf}{ \ddot {F} }
\newcommand{\wtG}{\wt{G}}
\newcommand{\dGammas}{\, d\Gamma_s}

\newcommand{\xu}{U}
\newcommand{\xz}{Z}
\newcommand{\xy}{Y}
\newcommand{\xl}{L}

\newcommand{\xR}{R}

\newcommand{\xX}{X}
\newcommand{\xY}{Y}
\newcommand{\xZ}{Z}
\newcommand{\xU}{U}
\newcommand{\xL}{L}

\newcommand{\cF}{{\mathcal F}}

\newcommand{\cG}{{\mathcal G}}

\newcommand{\FF}{{\mathbb F}}
\newcommand{\GG}{{\mathbb G}}
\newcommand{\RR}{{\mathbb R}}

\newcommand{\PP}{{\mathbb P}}
\newcommand{\EE}{{\mathbb E}}

\newcommand{\bigcdot}{\mathbin{{\hbox{\scalebox{.65}{$\bullet$}}}}}

\newcommand{\Keywords}[1]{\par\noindent{\small{\bf Keywords\/}: #1}}
\newcommand{\Class}[1]{\par\noindent{\small{\bf Mathematics Subjects Classification (2010)\/}: #1}}

\title{{{\Large \bf GENERALIZED BSDEs WITH RANDOM TIME HORIZON \\ IN A PROGRESSIVELY ENLARGED FILTRATION}} \vskip 40 pt }

\author{Anna Aksamit$\,^{a}$\footnote{The research of A. Aksamit was supported by the Australian Research Council DECRA Fellowship DE200100896.}
\,\,, Libo Li$\,^{b}$ and Marek Rutkowski$\,^{a,c}$\footnote{The research of M. Rutkowski was supported by the Australian Research Council Discovery Project DP200101550.} \\ \\
\\$^{a\,}$School of Mathematics and Statistics, University of Sydney \\ Sydney, NSW 2006, Australia \\ \\  
$^{b\,}$School of Mathematics and Statistics, University of New South Wales \\ Sydney, NSW 2052, Australia \\ \\ 
$^{c\,}$Faculty of Mathematics and Information Science, Warsaw University of Technology \\ 00-661 Warszawa, Poland \\ }

\date{\vskip 35 pt \today \vskip 30 pt}

\begin{document}
\maketitle
\begin{abstract}
We study generalized backward stochastic differential equations (BSDEs) up to a random time horizon $\xtheta$, which is not a stopping time,
under minimal assumptions regarding the properties of $\xtheta$.  In contrast to existing works in this area, we do not impose specific assumptions on the random time $\xtheta$ and we study the existence of solutions to BSDEs and reflected BSDEs with a random time horizon
through the method of reduction. In addition, we also examine BSDEs and reflected BSDEs with a l\`adl\`ag driver where the driver is allowed to have a finite number of common jumps with the martingale part.
\vskip 20 pt
\Keywords{BSDE, reflected BSDE, nonlinear evaluation, random time, enlargement of filtration}
\vskip 10 pt
\Class{60H99,$\,$91G99}
\end{abstract}

\newpage
\newpage

\section{Introduction} \label{sect1}

Our work is motivated by the arbitrage-free pricing of European and American style contracts with the counterparty credit risk
and, more generally, problems of mitigation of financial and insurance risks triggered by an extraneous event.
We focus on a study of BSDEs and reflected BSDEs up to a finite random time horizon $\xtheta$,
while imposing minimal assumptions on a random time $\xtheta$, which is not an $\FF$-stopping time with respect to a reference filtration $\FF$.
In contrast to existing papers (see, e.g., Ankirchner et al. \cite{ABE2010}, Kharroubi and Lim \cite{KL2014}, Cr\'epey and Song \cite{CS2015}, Dumitrescu et al. \cite{DQS2018}, and Grigorova et al. \cite{GQS}), we do not make any of simplifying assumptions frequently encountered in papers on the theory of progressive enlargement of filtration, such as: the immersion hypothesis, Jacod's equivalence hypothesis, the postulate {\bf (C)} of continuity of all $\FF$-martingales, or the postulate {\bf (A)} of avoidance of all $\FF$-stopping times by $\xtheta $. The only assumption we make is that the {\it Az\'ema supermartingale} of $\xtheta $ with respect to $\FF$
 (see Definition \ref{def2.2}) is a strictly positive process, although we also show that this assumption can be relaxed so our results apply to a larger class of random times, (see the class $\mathcal K$ in Section \ref{sect5.0}).

Stimulated by the recent paper by Choulli et al. \cite{CDV2020} on the martingale representation theorem in the {\it progressive enlargement} of a given filtration $\FF$ with observations of occurrence of a random time $\xtheta$, which is henceforth denoted as $\GG$, we study $\GG$-adapted BSDEs and $\GG$-adapted reflected BSDEs (RBSDEs) with an $\cF_\xtheta$-measurable terminal value at a random time horizon $\xtheta$. The crucial difference between the $\GG$ BSDEs and $\GG$ RBSDEs introduced in Definitions \ref{def3.1} and \ref{def3.4}, respectively, and various classes of BSDEs previously studied in the existing literature is that the driver is assumed to be a l\`agl\`ad process and the integrand against the compensated jump martingale $N^{o,\GG}$ given by equation \eqref{xeq2.1} is assumed to be $\FF$-optional, rather than $\GG$-predictable or, equivalently, $\FF$-predictable. 
Our main purpose in this work is to apply the {\it method of reduction} in the study of existence of a solution to the $\GG$ BSDE and $\GG$ RBSDE given by equations \eqref{eq3.2} and \eqref{eq3.5}, respectively. It should be acknowledged that the idea of reduction of a BSDE in a given filtration to a more tractable BSDE in a shrunken filtration has already been explored in papers by Kharroubi and Lim \cite{KL2014} and Cr\'epey and Song \cite{CS2015,CS2017} but the authors of these papers worked under various simplifying assumptions about the setup under study and did not consider the reflected case.

The first main contribution of this work is that we demonstrate that the idea of reduction can also be applied to the $\GG$ RBSDE \eqref{eq3.5}.
To be more precise, we show that the $\GG$ RBSDE with $\cF_\xtheta$-measurable terminal value can be solved up to a random horizon $\xtheta$ by first solving the corresponding reduced $\FF$ RBSDE and then constructing a solution to the original $\GG$ RBSDE by combining a solution to the $\FF$ RBSDE with an appropriate adjustment to the terminal value at time $\xtheta$. In particular, we analyze in detail the required adjustment at $\xtheta$ in the case when the driver of the $\GG$ RBSDE is a discontinuous l\`adl\`ag process. Furthermore, since the simplifying conditions {\bf (C)} and {\bf (A)} of continuity of all $\FF$-martingales and the avoidance of all $\FF$-stopping times
by $\xtheta$ are not imposed, we allow the reference filtration $\FF$ to support discontinuous martingales and, in addition, we also cover the situation where a random time $\xtheta$ has an overlap with $\FF$-stopping times. As a consequence, in the c\`adl\`ag case, the reduced $\FF$ BSDE obtained in our setup is such that the driver and the martingale part may share common jumps, which means that the reduced $\FF$ BSDE has a fairly general form that was not well studied in the existing literature.

Hence, as a second contribution of this work, we show how to construct a solution to the reduced $\FF$ BSDE for which the driver and the martingales appearing in a BSDE may have a finite number of common jumps. Our approach relies on a detailed analysis of the appropriate intermediate BSDE with a l\`agl\`ad driver. To illustrate our method more explicitly, we also show that a solution to an intermediate BSDE with a l\`agl\`ad driver can be obtained by adapting the existing results from Essaky et al. \cite{EHO2015} and Ren and El Otmani \cite{RE2010} who studied a particular class of BSDEs with a continuous driver.

We stress that the method presented here is not to solve directly $\GG$ BSDEs (or $\GG$ RBSDEs) through fixed point arguments by making appropriate technical assumptions on the solution space, the generator and the driver. Instead, our aim is to show that one can reduce the $\GG$ BSDE to a more manageable $\FF$ BSDE, which can be solved by making use of the solution to an intermediate BSDEs with l\`agl\`ad driver. Then we show that the solution of the intermediate BSDE with l\`agl\`ad driver can be obtained by a careful analysis of jumps, in a similar spirit to that of Essaky et al. \cite{EHO2015}, Klimsiak et al. \cite{KRS} and Confortola et al. \cite{CFJ2016}, and making use of existing results on BSDEs and RBSDEs with continuous drivers (such as those in \cite{EFO2017, EHO2015, FO2019,F2018,P1997,RE2010}) to deal with solutions on stochastic intervals between jumps.

The structure of the paper is as follows. In Section \ref{sect2}, we first introduce the notation and recall some auxiliary results from the theory of enlargement of filtration (see, e.g., Aksamit and Jeanblanc \cite{AJ2017}). Next, in view of recent works on RBSDEs with irregular barriers (see, e.g., Grigorova et al. \cite{GIOOQ2017,GIOQ2020} and Klimsiak et al. \cite{KRS}) and to demonstrate the generality of our methodology, we introduce in Definition \ref{def3.1} the notion of the l\`agl\`ad $\GG$ BSDE
(see also Definition \ref{def3.4} for the l\`agl\`ad $\GG$ RBSDE). In Section \ref{sect4}, we show how the reward process can be reduced and then discuss, in Sections \ref{sect5.0} and \ref{sect5}, some possible extensions of the setup given by Assumption \ref{ass3.1}.

Sections \ref{sect6}--\ref{sect8} are devoted to the issues of reduction of the $\GG$ BSDE \eqref{eq3.2} and construction of its solution.
We first demonstrate in Proposition \ref{pro4.1} that the $\GG$ BSDE can be effectively reduced to a corresponding equation in the filtration $\FF$. Subsequently, we then show in Proposition \ref{pro7.1} that a solution to the $\GG$ BSDE \eqref{eq3.2}
can be constructed from the corresponding reduced $\FF$ BSDE by first solving the coupled equations  \eqref{eq5.6}--\eqref{eq5.7}.
To examine the existence of a solution to the reduced $\FF$ BSDE, we first show that the $\FF$ BSDE \eqref{eq39.0}--\eqref{eq39.1} can be transformed
to the $\FF$ BSDE \eqref{eq40.0}--\eqref{eq40.1}, which is more tractable and whose solution can be used to address the issue of well-posedness of the coupled
equations \eqref{eq39.0}--\eqref{eq39.1}. Concrete situations where the coupled equations \eqref{eq39.0}--\eqref{eq39.1} possess a unique solution are studied in Section \ref{sect8} in the special case of a Brownian filtration (see Proposition \ref{pro9.1}).

In Sections \ref{sect9}--\ref{sect11}, we are concerned with analogous issues for $\GG$ RBSDEs and we first show that the method of reduction
can be used to reduce the $\GG$ RBSDE to the $\FF$ RBSDE. As shown in Proposition \ref{pro6.1},  the main new feature in the reflected case is that
the $\GG$-predictable reflection can be uniquely reduced to the $\FF$-predictable reflection, which is required to meet the appropriately modified Skorokhod conditions.
We then show in Proposition \ref{pro7.1} that, in principle, a solution to the $\GG$ RBSDE can be constructed from a solution to the reduced $\FF$ RBSDE.
The existence of a solution to the $\GG$ RBSDE in the Brownian case is studied in Section \ref{sect11} where we give in Proposition \ref{pro10.1}
sufficient conditions for the existence of a solution to the $\FF$ RBSDE in the case of the Brownian filtration $\FF$.

In Sections \ref{sect12}--\ref{sect13}, we deviate from the previous sections and, for given a filtration $\FF$, we focus on the BSDE \eqref{eq9.0} and the RBSDE \eqref{eq10.0}
with the feature that the driver is l\`agl\`ad and may shares jumps with the driving martingale. As mentioned before, even when the driver is c\`adl\`ag, there appears to be a gap in the existing literature on BSDEs when the driver shares common jumps with the driving martingale and thus we develop a jump-adapted method to solve BSDEs of this general form. Our method hinges on two steps. We first show, through a careful analysis of right-hand jumps, that the problem of solving the BSDE \eqref{eq9.0} on the whole interval $\llb 0, \tau \rrb$ can be addressed by solving a recursive system of c\`adl\`ag BSDEs \eqref{eq9.01} and then stitching together the solutions to that system. In the second step, we show in Proposition \ref{pro9.2} that a solution to a c\`adl\`ag BSDE can be obtained from a solution of an intermediate l\`agl\`ad BSDE \eqref{eq9.6}, which in turn can be solved by first solving a recursive system of c\`adl\`ag BSDEs \eqref{eq9.8} with a continuous driver, which are given on intervals defined by the right-hand jumps of l\`agl\`ad BSDE and, once again, appropriately aggregating these solutions. We argue that a reduction to the case of a continuous driver is important since it allow us to use existing results on well-posedness of BSDEs with a continuous driver. Concrete instances of our approach in a Brownian-Poisson filtration are presented in Examples \ref{ex9.2} and \ref{ex13.2}.

In the final section, we show that an analogous method can be applied to study the existence of a solution to the RBSDE \eqref{eq10.0} with a  l\`agl\`ad
driver using results for RBSDEs with a continuous driver. The main difference here is that we need to take care of the adjustment to the reflection process at the right-hand jumps and provide a rigorous check that the appropriate Skorokhod conditions are satisfied. The main result, Proposition \ref{pro10.2}, is complemented by an explicit illustration in a Poisson filtration given in Example \ref{ex10.1}.

\section{Setup and notation} \label{sect2}

Regarding the background knowledge, for the general theory of stochastic processes, we refer to He et al. \cite{HWY1992} and the reader interested in stochastic calculus for optional semimartingales is referred to Gal'\v{c}uk \cite{G1982}. For more details on the theory of random times and enlargement of filtration with applications to problems arising in financial mathematics (such as credit risk modeling or insider trading), the interested reader may consult the monograph by Aksamit and Jeanblanc \cite{AJ2017} and the recent paper by Jeanblanc and~Li~\cite{JL2020}.

We start by introducing the notation and recalling some fundamental concepts associated with modeling of a random time and the
associated notion of the progressive enlargement of a reference filtration.

We assume that a strictly positive and finite random time $\xtheta$, which is defined on a probability space $(\Omega, \cG, \PP)$, as well as some {\it reference filtration} $\FF$ are given. Then the enlarged filtration $\GG$ is defined as the
{\it progressive enlargement} of $\FF$ by observations of $\xtheta$ (see, e.g., \cite{AJ2017}) and thus a random time $\xtheta$, which is not necessarily
an $\FF$-stopping time and belongs to the set of all finite $\GG$-stopping times, denoted as $\cTg$.
We emphasize that the filtrations $\FF$ and $\GG$ are henceforth supposed to satisfy the usual conditions of $\PP$-completeness and right-continuity.

We will use the following notation for classes of processes adapted to the filtration $\FF$: \\

\vskip-10pt
\noindent  $\bullet$ $\cOff$, $\cPff$, $\overline{\mathcal{P}}(\FF)$ and $\cPrff$ are the classes of all real-valued, $\FF$-optional, $\FF$-predictable, $\FF$-strongly predictable and $\FF$-progressively measurable processes, respectively;\\
\noindent  $\bullet$ $\cOffd$, $\cPffd$, $\overline{\mathcal{P}}_d(\FF)$ and $\cPrffd$ are the classes of all $\RR^d$-valued, $\FF$-optional, $\FF$-predictable, $\FF$-strongly predictable and $\FF$-progressively measurable processes, respectively;\\
\noindent  $\bullet$ $\Martf$ (respectively, $\Martfl$) is the class of all $\FF$-martingales (respectively, $\FF$-local martingales);\\
\noindent  $\bullet$ $\Martft$ (respectively, $\Martflt$) is the class of all $\FF$-martingales (respectively, $\FF$-local martingales), which are stopped at $\xtheta $.

A stochastic process $X$ with sample paths possessing right-hand limits is said to be {\it $\FF$-strongly predictable}
if it is $\FF$-predictable and the process $X_+$ is $\FF$-optional (Definition 1.1 in \cite{G1982}).  
An analogous notation is used for various classes of $\GG$-adapted processes.
For instance, $\cPgg$ denotes the class of all $\GG$-predictable processes, $\Martglt$ is the class of all $\GG$-local martingales, which are stopped at the random time $\xtheta $, etc.  

In order to simplify the notation, we denote by $X\bigcdot Y$ the usual It\^o stochastic integral of $X$ with respect to a
(c\`adl\`ag) semimartingale $Y$, that is, $(X \bigcdot Y)_t := \int_{\rrb 0,t\rrb} X_s\,dY_s$, while we also write $(X \star Y)_t := \int_{\llb 0,t\llb} X_s\,dY_s$
so that the process $X \star Y$ is left-continuous as the integration is done over the interval $\llb 0,t\llb $. Due to the potential presence of a jump of $Y$ at time zero, we have that $(X\star Y)_t  = (X\bigcdot Y)_{t-} + X_0\Delta Y_0$  where, by the usual convention, $Y_{0-}=0$ so that $\Delta Y_0=Y_0$.

Let us recall from Gal'\v{c}uk \cite{G1982}  the notation pertaining to a pathwise decomposition of a l\`adl\`ag process. If $C$ is an $\FF$-adapted,
l\`adl\`ag process, then we write $C=C^c+C^d+C^g$ where the process $C^c$ is continuous, the c\`adl\`ag process $C^d$
equals $C^d_t:=\sum_{0\leq s\leq t}(C_{s}-C_{s-})$ and the c\`agl\`ad process $C^g$ is given by $C^g_t:=\sum_{0\leq s<t}(C_{s+}-C_s)$.
This also means that $C=C^r+C^g$ where the c\`adl\`ag process $C^r$ satisfies $C^r=C-C^g =C^c+C^d$. Notice that if $C$ is a c\`agl\`ad process,
then manifestly $C^d=0$ and thus $C^r=C^c$ is a continuous process. Similarly, if $C$ is a c\`adl\`ag process, then $C^g=0$ and
thus $C=C^r$.  For the sake of convenience, we denote by $C^g_+$ the c\`adl\`ag version of the c\`agl\`ad process $C^g$.


For a fixed random time $\xtheta$, we define the {\it indicator process} $A \in \cOgg$ by $A := \I_{\llb\xtheta,\infty\llb}$ so that $A_t=\I_\seq{\xtheta \leq t}$ for all $t\in \RR_+$ and we denote by $A^p$ (respectively, $A^o$) the dual $\FF$-predictable projection (respectively, the dual $\FF$-optional projection) of $A$. The BMO $\FF$-martingales $m$ and $n$ associated with the processes $A^o$ and $A^p$, respectively, are given by the following definition.

\bd \label{def2.1}
Let $m_t:=\EE(A^o_{\infty}\,|\,\cF_t)$ so that $m_{\infty}=A^o_{\infty}$ and let $n_t:=\EE(A^p_{\infty}\,|\,\cF_t)$ so that $n_{\infty}=A^p_{\infty}$.
\ed

\newpage

Following Az\'ema \cite{A1972}, we introduce the $\FF$-supermartingales $G$ and $\wtG$ associated with $\xtheta$.

\bd  \label{def2.2}
The c\`adl\`ag process $G\in \cOff $, which is given by the equality $G_t:=\PP (\xtheta > t\,|\,\cF_t)$, is called the {\it Az\'ema supermartingale} of $\xtheta $
with respect to $\FF$.  The l\`adl\`ag process $\wtG \in \cOff$, which is given by the equality $\wtG_t:=\PP ( \xtheta \geq t\,|\,\cF_t)$, is called the {\it
Az\'ema optional supermartingale} of $\xtheta$ with respect to~$\FF$.
\ed

For the reader's convenience, we recall some important properties of  Az\'ema supermartingales $G$ and $\wtG$ (see, e.g., \cite{AJ2017}).

\bl \label{lem2.1}
(i) We have that $G=n-A^p=m-A^o$ and $\wtG=m-A^o_{-}$ and thus
\[
G_t=\EE (A^p_{\infty }-A^p_{t}\,|\,\cF_t)=\EE (A^o_{\infty }-A^o_{t}\,|\,\cF_t),\quad \wtG_t=\EE (A^o_{\infty }-A^o_{t-}\,|\,\cF_t).
\]
(ii)  The processes $G$ and $\wtG$ satisfy $\wtG_-=G_-$ and $\wtG_+=G_+=G$.\\
(iii)  The inequality $\wtG\geq G$ holds and the equalities $\wtG-G={}^o(\I_{\llb \xtheta \rrb})=\Delta A^o$ and $\wtG-G_{-}=\Delta m$ are valid.
\el

Observe that the equality $G=n-A^p$ gives the Doob-Meyer decomposition in the filtration $\FF$ of the bounded $\FF$-supermartingale $G$.
From the classical theory of enlargement of filtration, it is well known that the $\GG$-martingale $N^{p,\GG}$ from the Doob-Meyer
decomposition in the filtration $\GG$ of the bounded $\GG$-submartingale $A$ can be represented as follows
\begin{equation} \label{eq2.2}
N^{p,\GG}:=A-\I_{\rrb 0,\xtheta\rrb}G^{-1}_{-}\bigcdot A^p .
\end{equation}

Furthermore, it was shown in Choulli et al. \cite{CDV2020} (see Theorem 2.3 therein) and Jeanblanc and Li \cite{JL2020} that the following process $N^{o,\GG}$ is a $\GG$-martingale with the integrable variation
\begin{equation} \label{xeq2.1}
N^{o,\GG}:=A-\I_{\rrb 0,\xtheta\rrb}\wtG^{-1}\bigcdot A^o=A-\I_{\rrb 0,\xtheta\rrb}\bigcdot\Gamma
\end{equation}
where the $\FF$ {\it hazard process} of $\xtheta$ is defined by $\Gamma := \wtG^{-1} \bigcdot A^o $ (see Jeanblanc and Li \cite{JL2020}).
The processes $N^{p,\GG}$ and $N^{o,\GG}$ are both known to belong to the class $\Martg$ but their properties are quite different; in particular, $N^{p,\GG}$ is not necessarily a {\it pure default martingale} (see Definition 2.2 in Choulli et al. \cite{CDV2020}) whereas $N^{o,\GG}$ has that property.

We will also make use of the following general result, which is due to Aksamit et al. \cite{ACJ2015}.

\bp \lab{pro2.1}
Let the $\FF$-stopping time $\wt \eta$ be defined by the following equality $\wt \eta:=\inf \{t\in \RR_+\,|\,\wtG_{t-}>\wtG_{t}=0\}$.
If $M$ is an $\FF$-local martingale, then
\begin{equation} \label{mardec}
M^\xtheta-\I_{\llb 0,\xtheta\rrb}\wtG^{-1}\bigcdot [M,m]+\I_{\llb 0,\xtheta\rrb}\bigcdot\big(\Delta M_{\wt \eta}\I_{[\![\wt \eta,\infty [\![}\big)^p
\end{equation}
is a $\GG$-local martingale stopped at $\xtheta$.
\ep

In particular, if $\wtG$ is a strictly positive process then for any $\FF$-local martingale $M$, we set
\begin{equation} \label{eq2.4}
\wt {M}:=M-\wtG^{-1}\bigcdot [M,m],
\end{equation}
so that the process $\wt M^\xtheta$ is a $\GG$-local martingale stopped at $\xtheta$. If, in addition, all $\FF$-martingales are continuous (that is, if the so-called Assumption {\bf (C)} holds, for instance, when $\FF$ is a Brownian filtration), then $\cOff =\cPff$ and thus the equalities $\wtG=G_-$ and $A^o=A^p$ are valid so that also $N^{o,\GG}=N^{p,\GG}$. Then equality \eqref{eq2.4} becomes
\begin{equation} \label{eq2.5}
\wt M=M-G^{-1}\bigcdot \langle M,n\rangle .
\end{equation}

\section{BSDEs with a random time horizon} \label{sect3}

Our study of BSDEs with a random time horizon is conducted within the following setup.

\bhyp \label{ass3.1}
We assume that we are given the following objects:\\
(i) a probability space $(\Omega,{\cal G},\PP )$ endowed with a filtration $\FF$; \\
(ii) a random time $\xtheta$ such that the Az\'ema supermartingale $G$ (hence also $G_{-}$ and the Az\'ema optional supermartingale $\wtG$) is a strictly positive process; \\
(iii) the class of all finite $\GG$-stopping times $\cTg$ where $\GG$ denotes the progressive enlargement of $\FF$ with a random time $\xtheta$; \\
(iv) the bounded processes $\xX, \xR \in \cOff$, which are used to define the bounded {\it reward process} $\Xgg\in\cOgg$ through the following expression
\begin{equation} \label{eq3.1}
\Xgg:=\xX\I_{\llb 0,\xtheta\llb}+\xR_\xtheta\I_{\llb\xtheta,\infty\llb};
\end{equation}
(v)  a real-valued $\GG$-martingale $N^{o,\GG}$ associated with $\xtheta$ and given by \eqref{xeq2.1}; \\
(vi) an $\RR^{d}$-valued, $\FF$-local martingale $M$, which is assumed to have the predictable representation property (PRP) for the filtration $\FF$;\\
(vii) an $\RR^{k}$-valued, $\GG$-adapted process $\wh D=(\wh D^1,\wh D^2,\dots,\wh D^k)$ where $\wh D^i$ is a linear combination of a l\`adl\`ag $\GG$-strongly predictable process of finite variation  and a l\`adl\`ag $\FF$-optional process of finite variation; \\
(viii) a mapping $\wh F = (\wh F^r , \wh F^g)$ where mappings $\wh F^r : \Omega \times \RR_+  \times \RR \times \RR^d \times \RR   \to \RR^k$
and $\wh F^g : \Omega \times \RR_+  \times \RR \times \RR^d \times \RR   \to \RR^k$
are such that, for any fixed $({y},{z},{u})\in \RR\times \RR^{d} \times \RR $, the process $(\wh F_t^r({y},{z},{u}))_{t\geq 0}$
belongs to $\mathcal{P}_k(\GG)$ and the process $(\wh F_t^g({y},{z},{u}))_{t\geq 0}$
belongs to $\mathcal{O}_k(\GG)$.
\ehyp

We are in a position to introduce a particular class of BSDEs with a random time horizon.
For the sake of brevity, they will be called $\GG$ BSDEs, as opposed to the associated $\FF$ BSDEs, which are introduced in Section \ref{sect6}.

\bd \label{def3.1}
For a fixed $\tgg \in \cTg$, we say that a triplet $(\Ygg,\Zgg,\Ugg)$ is a {\it solution} on $\llb 0, \tgg \wedge \xtheta \rrb$ to the $\GG$ BSDE
\begin{align} \label{eq3.2}
\Ygg_t&=\Xgg_{\tgg\wedge\xtheta}-\medint\int_{\,\rrb t,\tgg \wedge\xtheta\rrb}\wh F_s^r(\Ygg_s,\Zgg_s,\Ugg_s)\,d\wh D^r_s
-\medint\int_{\,\llb t,\tgg\wedge \xtheta\llb}\wh F_s^g(\Ygg_s,\Zgg_s,\Ugg_s)\,d\wh D^g_{s+} \nonumber \\
& \quad -\medint\int_{\,\rrb t,\tgg\wedge\xtheta\rrb}\Zgg_s\,d\wt{M}^\xtheta_s-\medint\int_{\,\rrb t,\tgg\wedge\xtheta\rrb}\Ugg_s\,d\Nogg_s
\end{align}
if $\Ygg\in\cOgg$ is a l\`agl\`ad process, the processes $\Zgg\in\cPggd$ and $\Ugg\in\cOff$
are such that the stochastic integrals in the right-hand side of \eqref{eq3.2} are
well defined and equality \eqref{eq3.2} is satisfied on the stochastic interval $\llb 0, \tgg \wedge \xtheta\rrb$.
\ed

The process $\wh{D}$ from Assumption \ref{ass3.1} (vii) and the mapping $\wh{F}= (\wh F^r, \wh F^g)$ from Assumption \ref{ass3.1} (viii) are henceforth called the {\it driver} and the {\it generator} of the $\GG$ BSDE, respectively. The processes $N^{o,\GG}$ and $\wt{M}^{\xtheta}$, given by equations \eqref{xeq2.1} and \eqref{eq2.4}, respectively, are orthogonal $\GG$-local martingales stopped at $\xtheta$ and they are referred to as {\it driving martingales.}
For explicit integrability conditions, which ensure that the stochastic integral $\Ugg \bigcdot \Nogg$ is a $\GG$-local martingale,
see Theorem 2.13 in Choulli et al. \cite{CDV2020}. To the best of our knowledge, the issue of well-posedness of the $\GG$ BSDE \eqref{eq3.2}
is not addressed in the existing comprehensive literature on BSDEs and thus our aim is to contribute to the theory of BSDEs by filling that gap.

%
%
In the next definition, we implicitly make the natural postulate of well-posedness of the $\GG$ BSDE \eqref{eq3.2} in a suitable
space of stochastic processes, which can be left unspecified at this stage.

\bd  \label{def3.2}
The {\it nonlinear evaluation} $\cEgg$ is the collection of mappings $\cEgg = \{\cEgg_{\sgg ,\tgg }\,|\, \sgg , \tgg \in \cTg ,\, \sgg \leq \tgg \}$
where for every $ \sgg , \tgg \in \cTg$ such that $\sgg \leq \tgg$ we have
$\cEgg_{\sgg,\tgg}(\Xgg_{\tgg}) := \Ygg_{\sgg\wedge \xtheta}$ where the triplet $(\Ygg,\Zgg,\Ugg )$
is a unique solution to the $\GG$ BSDE  \eqref{eq3.2} on the interval $\llb 0, \tgg\wedge\xtheta\rrb$.
\ed

\section{Reduction of the reward process}    \label{sect4}

Let $\cT$ denote the class of all finite $\FF$-stopping times and, for any fixed $\tau \in \cT$,  let the {\it stopped filtration} $\FF^\tau$ be given by $\FF^\tau := (\cF_{\tau \wedge t})_{t\geq 0}$. We will now examine the structure of the reward process $\Xgg$ specified by \eqref{eq3.1}. We claim that, for any $\tgg \in \cTg$, there exists $\tau \in \cT$ such that $\Xgg_{\tgg}=\Xgg_{\tau\wedge\xtheta}=\Xgg_\tau $. First, it is clear from \eqref{eq3.1} that the process $\Xgg$ is stopped at $\xtheta$ so that $\Xgg=\Xgg^\xtheta$, which immediately implies that $\Xgg_{\tgg}=\Xgg_{\tgg \wedge \xtheta }$. Hence, by using also the well known property that for any $\tgg \in \cTg$ there exists $\tau \in \cT$ such that $\tau \wedge \xtheta=\tgg\wedge \xtheta$, we obtain the following equalities
\begin{align*}
\Xgg_{\tgg}&=\Xgg_{\tgg\wedge\xtheta}=\xX_{\tgg\wedge\xtheta}\I_\seq{\tgg\wedge\xtheta<\xtheta}+\xR_\xtheta\I_\seq{\tgg\wedge\xtheta\geq\xtheta}= \xX_{\tau\wedge\xtheta}\I_\seq{\tau\wedge\xtheta<\xtheta}+\xR_\xtheta\I_\seq{\tau\wedge \xtheta\geq\xtheta}\\
&=\xX_{\tau}\I_\seq{\tau<\xtheta}+\xR_\xtheta\I_\seq{\tau\geq\xtheta}=\Xgg_{\tau\wedge\xtheta}=\Xgg_\tau
\end{align*}
so that $\Xgg_{\tgg}=\Xgg_{\tau\wedge\xtheta}=\Xgg_{\tau}$ for some stopping time $\tau \in \cT$, as was required to show.

\bl \label{lem4.3}
For any $\tau \in \cT $, there exists $\Xff(\tau)\in \mathcal{O}(\FF^\tau)$ such that the equality $\Xgg_\tau=\Xff_\xtheta(\tau )$ holds.
\el

\proof
It suffices to observe that
\begin{equation} \label{eq4.5}
\Xgg_{\tau}=\xR_\xtheta\I_{\llb \xtheta,\infty \rrb}(\tau)+\xX_{\tau}\I_{\rrb 0,\xtheta\llb}(\tau)
=\xR_\xtheta\I_{\llb 0,\tau\rrb}(\xtheta)+\xX_\tau\I_{\rrb\tau,\infty\llb}(\xtheta)=X_\xtheta(\tau)
\end{equation}
where, for any fixed $\tau $, the $\FF$-adapted process $\Xff(\tau)$ is given by
\begin{equation} \label{eq4.6}
\Xff(\tau ):=\xR\,\I_{\llb 0, \tau\rrb}+\xX_\tau\I_{\rrb\tau,\infty\llb}=\xR^\tau+(\xX_\tau-\xR_\tau)\I_{\rrb\tau,\infty\llb}.
\end{equation}
Since the processes $\xX$ and $\xR$ are assumed to be $\FF$-optional, by Lemma 3.53 in He et al. \cite{HWY1992}, the process $\Xff(\tau)$ is $\FF^\tau$-optional, although it is not a c\`adl\`ag process, in general.
\endproof

Since $\Xgg_{\tau}$ is $\cG_\tau$-measurable and $\Xff_\xtheta (\tau)$ is $\cG_{\xtheta}$-measurable, by part (3) of Theorem 3.4 in He et al. \cite{HWY1992}, the random variable $\Xgg_{\tau}=\Xff_\xtheta(\tau)$ is $\cG_{\tau\wedge\xtheta}$-measurable or, more precisely, it is $\cF^\tau_{\xtheta}$-measurable and $\cF^\tau_{\xtheta}\subset \cG_{\tau \wedge \xtheta}$. In view of equalities \eqref{eq4.5}, we will freely interchange $\Xgg_{\tgg},\,\Xgg_{\tau}$ and $\Xff_\xtheta(\tau)$.

\bp \label{pro3.1}
Given two $\GG$-stopping times $\sgg, \tgg$ such that $\sgg \leq \tgg$, there exists $\sigma \leq \tau$ where $\sigma , \tau \in \cT$ are such that $\sigma \wedge \xtheta=\sgg\wedge \xtheta,\, \tau \wedge \xtheta=\tgg\wedge \xtheta$ and $\cEgg_{\sgg, \tgg}(\Xgg_{\tgg \wedge \xtheta}) = \cEgg_{\sigma, \tau}(X_\xtheta(\tau))$.
\ep

\begin{proof}
First, we observe that there exists $\sigma$ such that $\sgg \wedge \xtheta = \sigma\wedge \xtheta$. Furthermore, if the inequality $\sigma \leq \tau $
fails to hold, then we can take $\sigma' = \sigma \wedge \tau$ and observe that $\sigma'\wedge \xtheta = \sigma\wedge \tau \wedge \xtheta  = \sgg\wedge \tgg \wedge \xtheta = \sgg\wedge \xtheta $. Using the fact that there exists $\sigma \leq \tau$ such that $\sigma \wedge \xtheta=\sgg\wedge \xtheta$ and $\tau \wedge \xtheta=\tgg\wedge \xtheta$, we obtain the following equalities
\begin{align*}
&\cEgg_{\sgg, \tgg}(\Xgg_{\tgg\wedge \xtheta})
= X_\xtheta(\tau) -\medint\int_{\,\rrb \sigma\wedge \xtheta, \tau\wedge \xtheta\rrb} \wh F_s^r(\Ygg_s,\Zgg_s,\Ugg_s)\,d\wh D^r_s \\
&-\medint\int_{\,\llb \sigma\wedge \xtheta,\tau \wedge \xtheta\llb} \wh F_s^g(\Ygg_s,\Zgg_s,\Ugg_s)\,d\wh D^g_{s+}
	 -\medint\int_{\,\rrb \sigma\wedge \xtheta,\tau \wedge \xtheta\rrb}\Zgg_s\,d\wt{M}^\xtheta_s-\medint\int_{\,\rrb \sigma\wedge \xtheta,\tau \wedge \xtheta\rrb}\Ugg_s\,d\Nogg_s
\end{align*}
and thus we conclude that $\cEgg_{\sgg, \tgg}(\Xgg_{\tgg\wedge \xtheta})= \cEgg_{\sigma, \tau}(X_\xtheta(\tau))$.
\end{proof}

\section{The case of a general Az\'ema supermartingale}     \label{sect5.0}


Let us make some comments on the possibility of relaxing Assumption \ref{ass3.1}(ii), that is, allowing the Az\'ema supermartingale $G$ of $\xtheta$ to hit zero. As an example, let us first consider a random time of the form $\xtheta = \xtheta'\wedge T_1$ where the Az\'ema supermartingale of $\xtheta'$ is a strictly positive process and $T_1$ is an $\mathbb{F}$-stopping time. Then the Az\'ema supermartingale of $\xtheta$ jumps to zero at $T_1$ and it is not hard to check that all results can be extended to that case by replacing the terminal time $\tau$ by $\tau' = \tau \wedge T_1$ and $\xtheta$ by $\xtheta'$. Since $G$ is a nonnegative supermartingale, we have that (see, e.g., Theorem 2.62 in He et al. \cite{HWY1992})
\[
\eta := 
\inf \{t\in \RR_+ \,|\, G_t = 0\} = \inf \{t\in \RR_+ \,|\, G_{t-} = 0\} = \lim_{n \to \infty} \eta_n
\]
where $\eta_n := \inf \{t \in \RR_+\,|\, G_t \leq 1/n \}$. It is known that $G=0$ on $\llb \eta , \infty \llb$ and, from Lemma 2.14 of \cite{AJ2017}, we have that $\llb 0, \xtheta \rrb \subset \{G_- >0\}$ and $\llb 0, \xtheta \llb \subset \{G >0\}= \llb 0, \eta \llb$ so that $\xtheta \leq \eta$.

In order to weaken Assumption \ref{ass3.1}(ii), we introduce the class of random times $\xtheta$ such that the Az\'ema supermartingale $G'$
of $\xtheta' := \xtheta_{\{\xtheta < \eta\}}$ is strictly positive on the interval $\llb 0, \eta \rrb$. Specifically, we set (recall that $\wtG_{-}=G_{-}$)
\begin{equation} \label{classK}
\mathcal{K} := \{ \xtheta \,|\, \wt G_{\eta}>0\}
\end{equation}
and we will argue that Assumption \ref{ass3.1}(ii) can be replaced by the postulate that $\xtheta$ belongs to $\mathcal{K}$.
It is clear that the class $\mathcal{K}$ is nonempty and contains not only all random times with a strictly positive
Az\'ema supermartingale, but also their minimum with any $\FF$-stopping time.
Observe that if a random time $\xtheta$ belongs to $\mathcal{K}$ then clearly
\[
\wt \eta:=\inf \{t\in \RR_+\,|\, \wtG_{t-}>\wtG_{t}=0\} = \eta_{\{G_{\eta-}>\wt G_\eta=0\}} = \infty
\]
and thus the term $\I_{\llb 0,\xtheta\rrb}\bigcdot\big(\Delta M_{\wt \eta}\I_{[\![\wt \eta,\infty [\![}\big)^p$ in the $\mathbb{G}$-semimartingale decomposition \eqref{mardec} of an arbitrary $\mathbb{F}$-local martingale $M$ is in fact null.

\bl \label{lemma5.1}
Let $\xtheta$ be a random time with the Az\'ema optional supermartingale $\wt G$ and let  $\xtheta' := \xtheta_{\{\xtheta < \eta\}}$.
Then  $\xtheta = \xtheta'\wedge \eta$ and the Az\'ema optional supermartingale $\wt G'$ of $\xtheta'$ is given by $\wt G' =\wt G^\eta +\I_{\rrb \eta ,\infty \llb}\bigcdot H$ where $H=\,^o(\I_{\{\xtheta=\eta\}})$.  Furthermore, we have that $(A')^o = (A^o)^\eta$ where $A':= \I_{\llb \xtheta', \infty\llb}$.
\el

\begin{proof}
We first observe that $\xtheta':=\xtheta_{\{\xtheta<\eta\}}=\xtheta\I_{\{\xtheta<\eta\}}+\infty\I_{\{\xtheta=\eta\}}$ where in the second equality we have used the inequality $\xtheta\leq\eta$. Using also the equalities $\{ \xtheta < \eta \}=\{ \xtheta' < \eta \}$
and $\{ \xtheta = \eta \}=\{ \xtheta' \geq \eta \}$, we obtain
\[
\xtheta = \xtheta \I_{\{\xtheta \leq \eta\}} = \xtheta \I_{\{\xtheta < \eta\}}+ \xtheta \I_{\{\xtheta = \eta\}}
= \xtheta' \I_{\{\xtheta < \eta\}}+ \eta \I_{\{\xtheta = \eta\}}
= \xtheta' \I_{\{\xtheta' < \eta\}}+ \eta \I_{\{\xtheta' \geq \eta\}} = \xtheta' \wedge \eta .
\]
Next, to compute the Az\'ema optional supermartingale of $\xtheta'$, we observe that for all $t> 0$,
\begin{align*}
\PP(\xtheta'<t \,|\,\cF_t)&=\PP(\xtheta<t,\xtheta<\eta\,|\,\cF_t)=\PP(\xtheta<t\,|\,\cF_t)-\PP(\xtheta<t,\xtheta=\eta\,|\,\cF_t)\\
&=1-\wt G_t\I_{\{t\leq\eta\}}-\PP(\xtheta=\eta\,|\,\cF_t)\I_{\{\eta<t\}}
\end{align*}
and hence $\wt G'=\wt G^\eta+\I_{\rrb \eta ,\infty \llb}\bigcdot H$ where $H=\,^o(\I_{\{\xtheta=\eta\}})$.
The last assertion follows from the uniqueness of the Doob-Meyer-Mertens decomposition of $\wt G'$.
It is worth noting that $(\wt G')^{\eta}=(\wt G)^\eta$ and $(G')^{\eta} = G\I_{\llb 0, \eta \llb} + \wt G_\eta \I_{\llb \eta, \infty\llb}$.
\end{proof}

It is important to observe that if $\xtheta \in \mathcal{K}$ then the supermartingale $G'$ and the optional supermartingale $\wt G'$ are strictly positive on $\llb 0,\tau'\rrb$ where $\tau':= \tau \wedge \eta$. Furthermore, the equalities $\Xgg_{\tgg}=\Xgg_{\tgg\wedge\xtheta}=\Xgg_{\tgg\wedge\xtheta'\wedge\eta}=X_{\xtheta'}(\tau')$ hold.
We conclude that, on the one hand, all our arguments used to address the case of a strictly positive Az\'ema supermartingale are still valid
when the pair $(\xtheta,\tau )$ is replaced by $(\xtheta',\tau')$. For instance, the BSDE in Proposition \ref{pro3.1} would not change since $\xtheta \leq \eta$, whereas in Section \ref{sect6} the terminal date can be changed to $\tau'=\tau\wedge\eta$. On the other hand, however, if a random time $\xtheta$ is not in $\mathcal{K}$ then technical issues involving either an explosion of integrals or ill-defined terminal condition at time $\eta$ may arise.

\section{Extended terminal condition}  \label{sect5}

Let us make some comments on a possible extension of the terminal condition in the $\GG$ BSDE. Since the processes $\xX$ and $\xR$ are assumed to be $\FF$-optional, Lemma \ref{lem4.3} implies that $X(\tau)$ is an $\FF^\tau$-optional process and $X_\xtheta(\tau)$ is $\cF^{\tau}_{\xtheta}$-measurable. This implies that in our formulation of the $\GG$ BSDE \eqref{eq3.2} and the $\GG$ RBSDE \eqref{eq3.5} the terminal condition is measurable with respect to $\cF^\tau_{\xtheta} \subset \cG_{\tau\wedge \xtheta}$. It is also worth noting that
\[
\sigma(V^\tau_\xtheta \,|\,  V \in \mathcal{O}(\FF)) = \cF_{\tau\wedge \xtheta} \subset \cF^{\tau}_{\xtheta} := \sigma(V_\xtheta \,|\,  V \in \mathcal{O}(\FF^\tau)).
\]
Since we do not consider all $\cG_{\tau\wedge \xtheta}$-measurable terminal conditions, the multiplicity in the martingale representation property established in Theorem 2.22 of Choulli et al. \cite{CDV2020} can be taken to be equal to two, which in fact gives a partial motivation for Definition \ref{def3.1} of a solution to the $\GG$ BSDE.

In general, the multiplicity in Theorem 2.22 of Choulli et al. \cite{CDV2020} is equal to three and thus it would be possible to consider, when a (bounded) terminal condition $\zeta$ is $\cG_{\tau\wedge \xtheta}$-measurable, a more general $\GG$ BSDE driven by $\wt M^{\xtheta},\, \Nogg$ and a pure jump martingale yielding an additional `correction term' at the terminal time $\tau\wedge \xtheta$. To be more specific, in view of Proposition \ref{pro3.1}, one could study the extended $\GG$ BSDE of the form
\begin{align}
\Ygg_t & = \zeta - \medint\int_{\,\rrb t,\tau\wedge \xtheta\rrb} \wh F^r_s(\Ygg_s,\Zgg_s, \Ugg_s, \whJ_s(\tau))\, d\wh D^r_s - \medint\int_{\,\llb t,\tau\wedge \xtheta \llb}   \wh F^g_s(\Ygg_s,\Zgg_s, \Ugg_s, \whJ_s (\tau) )\, d\wh D^g_{s+}  \nonumber \\
	   & \quad - \medint\int_{\,\rrb t,\tau\wedge \xtheta\rrb}  \Zgg_s \,d\wt M_s^\xtheta - \medint\int_{\,\rrb t,\tau \wedge \xtheta\rrb}  \Ugg_s \, d\Nogg_s - \medint\int_{\,\rrb t,\tau\wedge \xtheta \rrb}  \whJ_s(\tau)\, dA_s  \label{kbsde}
\end{align}
where $\Zgg$ is $\FF$-predictable, $\Ugg$ is $\FF$-optional, $\whJ(\tau)$ is $\FF^\tau$-progressively measurable and $\mathbb{E}(\whJ_{\xtheta}(\tau)\,|\,\cF_\xtheta) = 0$. It is not difficult to check that the last condition implies that $\whJ(\tau) \bigcdot A$ is a $\GG$-local martingale, provided that a suitable integrability condition is satisfied by $\whJ(\tau)$.

A detailed study of the $\GG$ BSDE given by equation \eqref{kbsde} is beyond the scope of this work since its practical applications are unclear.
Let us only point out that if the generator $\wh F$ in \eqref{kbsde} does not depend on $\whJ$, then one can formally reduce \eqref{kbsde} to \eqref{eq3.2} and show that a solution to \eqref{kbsde} can be obtained from a solution to \eqref{eq3.2}. To this end, we will need the following auxiliary result.

\bl \label{lem3.1}
Assume that $\zeta$ is bounded and $\cG_{\tau\wedge \xtheta}$-measurable. Then \\
(i) there exists a process $X'(\tau) \in \mathcal{P}r(\FF^\tau)$ such that $\zeta = X'_{\xtheta}(\tau)$,\\
(ii) there exists a process $X(\tau) \in \mathcal{O}(\FF^\tau)$ such that $X_{\xtheta}(\tau) = \EE(X'_{\xtheta}(\tau)\,|\,\cF_\xtheta)$.
\el

\begin{proof}
To show the first assertion, we note that since $\zeta$ is $\cG_{\tau\wedge \xtheta}$-measurable, there exists a process $\wh H \in \cOgg$ such that $\zeta = \wh H_{\tau\wedge \xtheta}$ and thus also, by Proposition 2.11 in Aksamit and Jeanblanc \cite{AJ2017}, a process $H\in \cOff$ such that $\wh H\I_{\llb 0, \xtheta \llb} = H\I_{\llb 0, \xtheta \llb}$.

 Furthermore, since $\cG_{\tau\wedge \xtheta} \subset \cG_{\xtheta} = \cF_{\xtheta+}$, there exists a process $H' \in \cPrff$ such that $\zeta = \wh H_{\tau \wedge \xtheta} = H'_{\xtheta}$ (see Lemma B.1 in Aksamit et al. \cite{ACDJ2017}, which is obtained by modifying Proposition 5.3 (b) in Jeulin \cite{J1980}). We thus have the equalities
 \[
 H'_{\xtheta}\I_{\{\tau < \xtheta\}}  = \wh H_{\tau} \I_{\{\tau < \xtheta\}} = H_{\tau} \I_{\{\tau < \xtheta\}}
 \]
 and we can define the $\FF^\tau$-progressively measurable process
\begin{equation} \label{eqX'}
X'(\tau) := H'\I_{\llb 0,\tau\rrb}+H_{\tau}\I_{\rrb \tau,\infty\llb}
\end{equation}
which satisfies $X'_{\xtheta}(\tau) = \zeta$. For the second assertion, we note that since $\I_{\llb 0, \tau \rrb}$ and $H_{\tau} \I_{\rrb \tau, \infty \llb}$ belong to $\cOff$, we have
\[
\EE(X'_\xtheta(\tau)\,|\,\cF_\xtheta) = \EE(H'_\xtheta\,|\,\cF_\xtheta)\I_{\llb 0,\tau \rrb}(\xtheta)+H_{\tau} \I_{\rrb \tau,\infty \llb}(\xtheta).
\]
By Proposition 2.21 in Choulli et al. \cite{CDV2020} there exists an $\FF$-optional process $X$ such that $X_\xtheta = \EE(H'_\xtheta\,|\,\cF_\xtheta)$.
It now suffices to set $X(\tau):=X\I_{\llb 0, \tau \rrb}+H_{\tau}\I_{\rrb \tau, \infty \llb}$ and observe that the process $X(\tau )$ belongs to $\mathcal{O}(\FF^\tau)$.
\end{proof}

By applying Lemma \ref{lem3.1} to an integrable, $\cG_{\tau \wedge \xtheta}$-measurable random variable $\zeta$, we can rewrite
$\zeta = X'_{\xtheta}(\tau) = X'_{\xtheta}(\tau)-X_{\xtheta}(\tau)+X_{\xtheta}(\tau)$. Therefore, in view of \eqref{eqX'}, we have
$\whJ(\tau):= X'(\tau)-X(\tau)=(H'-X)\I_{\llb 0,\tau\rrb}$,
which shows that $\whJ(\tau)=\whJ(\tau)\I_{\llb 0,\tau\rrb}$.

\newpage

 Then the BSDE \eqref{kbsde} becomes
\begin{align*}
\Ygg_t &= X_{\xtheta}(\tau) - \medint\int_{\,\rrb t,\tau\wedge \xtheta\rrb} \wh F^r_s(\Ygg_s,\Zgg_s, \Ugg_s)\,d\wh D^r_s
- \medint\int_{\,\llb t,\tau \wedge \xtheta\llb} \wh F^g_s(\Ygg_s,\Zgg_s, \Ugg_s)\, d\wh D^g_{s+}  \nonumber \\
&\quad - \medint\int_{\,\rrb t,\tau \wedge \xtheta\rrb}  \Zgg_s\, d\wt M_s^\xtheta - \medint\int_{\,\rrb t,\tau\wedge \xtheta\rrb}\Ugg_s\,d\Nogg_s
\end{align*}
where we have used the equalities
\[
\medint\int_{\,\rrb t,\tau\rrb} \whJ_s(\tau)\, dA_s = \whJ_{\xtheta}(\tau) \I_{\{t< \xtheta \leq \tau\}} = \whJ_{\xtheta}(\tau).
\]
We conclude that if $(\Ygg, \Zgg, \Ugg)$ is a solution to the $\GG$ BSDE \eqref{eq3.2}, but with $\tgg$ replaced by $\tau$, then $(\Ygg, \Zgg, \Ugg, \whJ(\tau))$ solves the BSDE \eqref{kbsde} with a $\cG_{\tau\wedge \xtheta}$-measurable terminal condition. Similar arguments can be applied to the case of the RBSDE \eqref{eq3.5}. However, in the case of generators depending on $\whJ(\tau)$, the proper form of the adjustment to the terminal condition would be more complicated and its computation would involve the generators $\wh F^r$ and $\wh F^g$.

\section{Reduction of a solution to the $\GG$ BSDE}   \label{sect6}

Our goal is to show that the BSDE \eqref{eq3.2} has a solution, which can be obtained in two steps. In the reduction step, the filtration is shrunk from $\GG$ to $\FF$ and the BSDE \eqref{eq3.2} is analyzed through an associated reduced BSDE in the filtration $\FF$. In the construction step, we show that a solution to the reduced BSDE can be lifted from $\FF$ to $\GG$ in order to obtain a solution to the BSDE \eqref{eq3.2}. Notice that in Sections \ref{sect6} and \ref{sect7} the random times $\xtheta $ and $\tgg$ are fixed throughout.

We first establish some preliminary lemmas related to the concept of shrinkage of filtration. In the main result of this section, Proposition \ref{pro4.1}, we give an explicit representation for the $\FF$ BSDE associated with the $\GG$ BSDE \eqref{eq3.2}. We work here under Assumption \ref{ass4.1}, which will be relaxed in Section \ref{sect7} where an explicit construction of a solution to the $\GG$ BSDE is proposed and analyzed.

\bhyp \label{ass4.1}
A solution $(\Ygg,\Zgg,\Ugg)$ to the BSDE \eqref{eq3.2} exists on the stochastic interval $\llb 0, \tgg\wedge\xtheta\rrb$ or,
equivalently, on the interval $\llb 0, \tau \wedge\xtheta\rrb$ where $\tau \in \cT$ is such that $\tau \wedge \xtheta=\tgg\wedge \xtheta$.
\ehyp

Our present goal is to analyze the consequences of Assumption \ref{ass4.1}. We start by recalling that there exist
a unique $\FF$-optional process ${Y}$ and a unique $\FF$-predictable process $\xZ$ such that the equalities $\Ygg\I_{\llb 0,\xtheta\llb}=\xY\I_{\llb 0,\xtheta\llb}$ and $\Zgg \I_{\llb 0,\xtheta\rrb}= \xZ\I_{\llb 0,\xtheta\rrb}$ are valid. Moreover, ${Y}_{\tau}=\xX_{\tau}$ and the process ${Y}$ and ${Z}$ are given by
\begin{equation} \label{eq4.1}
\xY={}^o\big(\I_{\llb 0,\xtheta\llb}\Ygg\big)G^{-1},\quad \xZ ={}^p\big(\I_{\llb 0,\xtheta\rrb}\Zgg\big)G^{-1}_-.
\end{equation}
Similarly, in view of Assumption \ref{ass3.1}(vii) and Lemma \ref{lem6.1} below, there exists an right continuous $\FF$-adapted process $D^r$ and a left-continuous $\FF$-adapted process $D^g$ such that $\wh D^r \I_{\llb 0,\xtheta \rrb} = D^r \I_{\llb 0,\xtheta \rrb}$ and $\wh D^g \I_{\llb 0,\xtheta \rrb} = D^g \I_{\llb 0,\xtheta \rrb}$.
Finally, it is clear that $\wh X\I_{\llb 0,\xtheta \llb} = X\I_{\llb 0,\xtheta \llb}$. We shall refer to $\tau$, $\xY$, $\xZ$, $D^r$, $D^g$ and $X$ as the $\FF$-{\it reduction} of $\tgg$, $\Ygg$, $\Zgg$, $\wh D^r$, $\wh D^g$ and $\wh X$.

In the following, we slightly abuse the notation and we again denote by $\xY$ and $\xZ$ the stopped processes $\xY:=\xY^{\tau}$ and $\xZ:=\xZ^\tau $. Recall that the component $\Ugg$ in a solution to \eqref{eq3.2} is assumed to be an $\FF$-optional process and thus $\Ugg$ is
equal to its $\FF$-reduction $\xU$ so that, trivially, $\Ugg =\xU$ and, once again, we will write $U := \xU^\tau$.

To show more explicitly how the process $Y$ is computed, we observe that
\[
\cEgg_{\cdot, \tgg}(\Xgg_{\tgg})\I_{\llb 0, \tau \rrb}\I_{\llb 0, \xtheta \llb}
 = \cEgg_{\cdot, \tgg}(\Xgg_{\tgg})\I_{\llb 0, \tgg \rrb}\I_{\llb 0, \xtheta \llb}
\] and, in view of Proposition \ref{pro3.1}, there exists $Y \in \cOff$ such that we have the equalities
\[
\cEgg_{\cdot, \tgg}(\Xgg_{\tgg})\I_{\llb 0, \tau \rrb}\I_{\llb 0, \xtheta \llb}  = \cEgg_{\cdot, \tau}(X_\xtheta(\tau))\I_{\llb 0, \tau \rrb}\I_{\llb 0, \xtheta \llb} = Y\I_{\llb 0, \tau \rrb}\I_{\llb 0, \xtheta \llb}.
\]

\newpage

Therefore, by applying the $\FF$-optional projection operator, we obtain
\[
^o\big(\cEgg_{\cdot,\tgg}(\Xgg_{\tgg})\I_{\llb 0,\tgg\rrb}\I_{\llb 0, \xtheta \llb}\big)\I_{\llb 0, \tau \rrb}
=\,^o\big(\cEgg_{\cdot,\tau}(X_\xtheta(\tau))\I_{\llb 0,\tau \rrb}\I_{\llb 0,\xtheta \llb}\big)\I_{\llb 0, \tau \rrb}
= YG\I_{\llb 0, \tau \rrb}.
\]

A general representation of $Y$ can then be obtained on $\llb 0, \tau \rrb$ by noticing that for any $\FF$-stopping time $\sigma$
\[
Y_\sigma G_\sigma \I_{\{ \sigma \leq  \tau \}} =
\EE \big(\cEgg_{\sigma,\tau}(X_\xtheta(\tau))\I_{\{\sigma \leq \tau\}}\I_{\{\sigma<\xtheta \}}\,|\, \cF_\sigma\big)\I_{\{ \sigma \leq  \tau \}}.
\]
Our next goal is provide a more explicit computation of the right-hand side in the above equality (see Lemma \ref{lem4.4}).

\brem
Suppose that Assumption \ref{ass3.1}(ii) is relaxed and we postulate instead that $\xtheta \in \mathcal{K}$ where the class $\mathcal{K}$ is defined by \eqref{classK}. Then the modified terminal condition would be $\Xgg_{\tgg\wedge\xtheta'\wedge \eta} = \Xgg_{\tau\wedge\xtheta'\wedge \eta}$ and the reduced terminal condition would become $X_{\tau \wedge \eta}G'_{\tau \wedge \eta}$. Finally, the terminal condition for $Y$ would be $X_{\tau \wedge \eta}$ instead of $X_\tau$. Hence, it would be enough to replace $\tau$ with $\tau' = \tau \wedge \eta$ and study the $\FF$-BSDE on the interval $\llb 0, \tau'\rrb$, rather than $\llb 0, \tau\rrb$.
\erem

The following result can be deduced from Proposition 2.11 in Aksamit and Jeanblanc \cite{AJ2017}.

\bl \label{lem4.2a}
For every $(y,z,u)\in \RR\times \RR^{d}\times \RR $ there exists an $\RR^{k}$-valued, $\FF$-predictable process $F^r(y,z,u)$ such that $\wh F^r_t(y,z,u)\I_\seq{\xtheta\geq t}=F^r_t(y,z,u)\I_\seq{\xtheta\geq t}$ for every $t\geq 0$. For every $(y,z,u)\in \RR\times \RR^{d}\times \RR $ there exists an $\RR^{k}$-valued, $\FF$-optional process $F^g(y,z,u)$ such that $\wh F^g_t(y,z,u)\I_\seq{\xtheta > t}=F^g_t(y,z,u)\I_\seq{\xtheta >  t}$ for every $t\geq 0$.
\el

To reduce the driver $\wh D$ and later the reflection in the $\GG$ RBSDE, we prove the following result. Notice that a similar result was established in Jeanblanc et al. \cite{JLS2018} in the case where the partition of the space $\Omega \times [0,\infty[$ was independent of time.

\bl \label{lem6.1}
Let $\wh D = \wh D^r + \wh D^g$  be an $\GG$-adapted, l\`agl\`ad, increasing process. Then there exists an $\FF$-optional, c\`adl\`ag, increasing process $D^r$ and an $\FF$-predictable, c\`agl\`ad, increasing process $D^g$ such that $D^r = \wh D^r$ on $\llb 0,\xtheta \llb$ and $D^g = \wh D^g$ on $\llb 0,\xtheta \rrb$. If $\wh D$ is a $\GG$-strongly predictable increasing process, then $D^r$ can be chosen such that it is an
$\FF$-predictable, c\`adl\`ag, increasing process and $D^r = \wh D^r$ on $\llb 0,\xtheta\rrb$.
\el

\begin{proof}
Since $\wh D^r$ belongs to the class $\cOgg$, there exists an $\FF$-optional process $D^r$ such that $\wh D^r\I_{\llb 0, \xtheta \llb}=D^r\I_{\llb 0,\xtheta \llb}$ (see the first equality in \eqref{eq4.1}). Since the optional projection of a c\`adl\`ag processes is again a c\`adl\`ag process, the process $D^r$ is c\`adl\`ag on the set $\{G>0\}=\Omega \times [0,\infty[$ where the last equality is clear since we have assumed that $G$ is strictly positive.

To show that the process $D^r$ is increasing, we observe that, for every $s\leq t$,
\[
D^r_t\I_{\{\xtheta > t\}}= \wh D^r_t\I_{\{\xtheta > t\}}
\geq \wh D_s^r\I_{\{\xtheta > t\}}= \wh D_s^r\I_{\{\xtheta > s\}}\I_{\{\xtheta > t\}}= D_s^r\I_{\{\xtheta > t\}}.
\]
Then, by taking the $\cF_{t}$ conditional expectation of both sides, we deduce that the process $D^r$ is increasing on the set $\{G>0\}=\Omega \times [0,\infty[$.

Furthermore, since the process $\wh D^g$ is c\`agl\`ad and thus belongs to the class $\cPgg$, there exists an $\FF$-predictable process $D^g$ such that $\wh D^g\I_{\llb 0, \xtheta \rrb}=D^g\I_{\llb 0,\xtheta \rrb}$ (see the second equality in \eqref{eq4.1}). The rest of the proof is similar to the case of $D^r$ except that we now use the properties of the $\FF$-predictable projection, rather than the $\FF$-optional projection.

Finally, in the case where $\wh D$ is $\FF$-strongly predictable, from the decomposition $\wh D = \wh D^r + \wh D^g$ and the fact that $\wh D$ and $\wh D^g$ belong to $\cPgg$, we deduce that $\wh D^r$ belongs to $\cPgg$. Thus there exists an $\FF$-predictable process $D^r$ such that $\wh D^r\I_{\llb 0, \xtheta \rrb}=D^r\I_{\llb 0,\xtheta \rrb}$. This implies that $\wh D^r\I_{\llb 0, \xtheta \llb}=D^r\I_{\llb 0,\xtheta \llb}$ and, by taking the $\FF$-optional projection, we deduce from similar arguments as before, that on the set $\{G>0\}$ the process $D^r$ is increasing and c\`adl\`ag.
\end{proof}

\brem
Clearly if the process $\wh D$ is $\FF$-adapted then the equality $\wh D^r = D^r$ holds everywhere and not only on $\llb 0,\xtheta\llb$. We remark here that an $\FF$-adapted driver $\wh D$ can have certain practical interpretations. For example, one can take $\wh D$ to be the hazard process, that is $\wh D = \Gamma = \wt G^{-1}\bigcdot A^o$, and this can be interpreted as a way to introduce ambiguity in the recovery and the default intensity (see, for instance, Fadina and Schmidt \cite{FS2019}).
\erem


The next result is an immediate consequence of Lemma \ref{lem6.1} and equations \eqref{eq3.1} and \eqref{eq3.2}.
To alleviate the notation, we will frequently write $\wh F_s^r(\cdot)=\wh F_s^r(\cdot ,\Zgg_s,\Ugg_s),$
$\wh F_g^r(\cdot)=\wh F_s^g(\cdot,\Zgg_s,\Ugg_s), F_s^r(\cdot )=F_s^r(\cdot,\xZ_s,\xU_s)$ and $F_s^g(\cdot)=  F_s^g(\cdot,\xZ_s,\xU_s)$.

\bl \label{lem4.2}
The following equalities are satisfied, for every $t\in \RR_+$ on the event $\seq{t\leq\tau}\cap \seq{t<\xtheta}$,
\begin{align} \label{eq4.3}
&\EE \Big( \I_\seq{\xtheta>t}\medint\int_{\,\rrb t,\tau\wedge\xtheta\rrb}\wh F_s^r(\Ygg_s)\,d\wh D_s^r\,\big|\,\cF_t\Big)
=\EE \Big(\medint\int_{\,\rrb t,\tau\rrb} F_s^r(\xY_s)\,dD_s^r\,\big|\,\cF_t\Big) \nonumber \\
&\qquad + \EE \Big( \medint\int_{\,\rrb t,\tau\rrb} \big(F_s^r(\xR_s)-F_s^r(\xY_s)\big) \Delta D_s^r dA^o_s \,\big|\,\cF_t \Big)
\end{align}
and, on the event $\seq{t < \tau}\cap \seq{t<\xtheta}$,
\[
\EE \Big(\I_\seq{\xtheta>t}\medint\int_{\,\llb t,\tau\llb}\I_\seq{\xtheta > s}\wh F_s^g(\Ygg_s)\,d\wh D^g_{s+}\,\big|\,\cF_t\Big)
= \EE \Big(\medint\int_{\,\llb t,\tau\llb} G_s F_s^g(\xY_s)\,dD^g_{s+}\,\big|\,\cF_t\Big).
\]
\el

\proof
Using Lemma \ref{lem6.1}, the equalities $\Ygg_{\xtheta}\I_\seq{\xtheta \leq \tau}=R_{\xtheta}\I_\seq{\xtheta \leq \tau}$ and $\PP (\xtheta=s\,|\,\cF_s)=\wtG_s-G_s=\Delta A^o_s$ and noticing that $\seq{\xtheta \geq s}\subset \seq{\xtheta > t}$ for $s>t$, we obtain
\begin{align*}
&\EE \Big(\I_\seq{\xtheta>t}\medint\int_{\,\rrb t,\tau\rrb}\I_\seq{\xtheta \geq s}\wh F_s^r(\Ygg_s)\,d\wh D^r_s\,\big|\,\cF_t\Big)\\
&=\EE \Big(\medint\int_{\,\rrb t ,\tau\rrb}\I_\seq{\xtheta>s} F_s^r (\xY_s)\,dD^r_s\,\big|\,\cF_t\Big)
+\EE\Big(\medint\int_{\,\rrb t,\tau\rrb}\I_\seq{\xtheta =s} F_s^r (\xR_s)\,dD^r_s\,\big|\,\cF_t\Big) \\
&=\EE \Big(\medint\int_{\,\rrb t,\tau\rrb}G_s F_s^r(\xY_s)\, dD^r_s\,\big|\,\cF_t\Big)+
\EE \Big(\medint\int_{\,\rrb t,\tau\rrb} F_s^r (\xR_s)\Delta A^o_s\, dD^r_s\,\big|\,\cF_t\Big)\\
&=\EE \Big(\medint\int_{\,\rrb t,\tau\rrb} \wt G_s F_s^r(\xY_s)\, dD^r_s\,\big|\,\cF_t\Big)
+ \EE \Big(\medint\int_{\,\rrb t,\tau\rrb}(F_s^r(\xY_s)-F_s^r(\xR_s))\Delta D_s^r\,dA^o_s \,\big|\,\cF_t\Big).
\end{align*}

Similarly, again from Lemma \ref{lem6.1}, we have
\[
\EE \Big(\I_\seq{\xtheta>t}\medint\int_{\,\llb t,\tau\llb}\I_\seq{\xtheta > s}\wh F_s^g(\Ygg_s)\,d\wh D^g_{s+}\,\big|\,\cF_t\Big)
= \EE \Big(\medint\int_{\,\llb t,\tau\llb} G_s F_s^g(\xY_s)\,dD^g_{s+}\,\big|\,\cF_t\Big),
\]
which gives the required result.
\endproof

To simplify further computations we define, for every $(y,z,u)\in \RR\times \RR^{d}\times \RR $ and $t\geq 0$,
\begin{equation} \label{eq4.4}
\wtf^r_t(y,z,u):= F^r_t(y,z,u)+ \big(F^r_t(R_t,z,u)-F^r_t(y,z,u)\big)G^{-1}_t\Delta A^o_t .
\end{equation}
By combining Lemmas  \ref{lem4.3} and \ref{lem4.2} with equality \eqref{eq3.2}, we obtain the following result.

\bl \label{lem4.4}
The process $(Y,Z,\xU)$ satisfies on $\llb 0,\tau\rrb$
\[
Y_t=G^{-1}_t\,\EE\Big(\xX_{\tau}G_{\tau}+\medint\int_{\,\rrb t,\tau \rrb}\xR_s\,dA^o_s -\medint\int_{\,\rrb t,\tau \rrb} \wt G_s \wtf_s^r \,dD^r_s -\medint\int_{\,\llb t,\tau \llb} G_s F_s^g\,dD^g_{s+}\,\,\big|\,\cF_t\Big)
\]
where we denote $\wtf_s^r =\wtf_s^r (\xY_s,\xZ_s,\xU_s)$ and $F_s^g=F_s^g(\xY_s,\xZ_s,\xU_s)$.
\el

\begin{proof}
For any fixed $\FF$-stopping time $\tau$, we denote by $\xNX(\tau)$ the $\FF$-martingale given by
\begin{equation}  \label{eq4.8}
\xNX_t(\tau):=\EE \big(\xX_{\tau}G_{\tau}+ (\xR \bigcdot A^o)_\tau
-(\wt G \wtf^r \bigcdot D^r)_\tau - (G F^g \star D^g_+)_\tau \,|\,\cF_t\big).
\end{equation}
Then the $\FF$-optional process $Y$ has the following representation on $\llb 0,\tau\rrb$
\begin{equation} \label{eq4.9}
Y_t=G^{-1}_t\,\big(\xNX_t(\tau)-(\xR \bigcdot A^o)_t +(\wt G\wtf^r\bigcdot D^r)_t + (G F^g\star D^g_+)_t  \big)
\end{equation}
with $Y_{\tau}=\xX_{\tau }$ and thus the asserted equality holds.
\end{proof}

For brevity, we set $C:= \wt G\wtf^r \bigcdot D^r + GF^g\star D^g_{+}$ and we note that equality \eqref{eq4.9} can be rewritten as follows
\begin{equation}  \label{eq4.10}
Y=G^{-1}\big(\xN(\tau ) - \xR \bigcdot A^o +C\big).
\end{equation}
In addition, we define $\wt{m}:=m -\wtG^{-1}\bigcdot [m,m]$ and
\[
\wt{\xN}(\tau):=\xN(\tau )-\wtG^{-1}\bigcdot  [\xN(\tau) ,m].
\]
To express the dynamics of the process $Y$ in terms of $\wt{m}$ and $\wt{\xN}(\tau )$, we will use the following immediate consequence
of Lemma \ref{lem11.3} from the appendix.

\bl \label{lem4.5}
If $C=C^r+C^g$ is a l\`agl\`ad process of finite variation and the process $Y$ is given by
$Y=G^{-1}(K- \xR \bigcdot A^o +C)$ for some $\FF$-martingale $K$, then
\[
\xY=\xY_0+ \wtG^{-1} \bigcdot C^r+G^{-1}\star C^g_{+}- (\xR-\xY)\bigcdot \Gamma -\xY_{-}G_{-}^{-1}\bigcdot \wt m+G^{-1}_{-}\bigcdot \wt K
\]
where $\wt K:=K-\wtG^{-1}\bigcdot  [K,m]$.
\el

By applying Lemma \ref{lem4.5} to equality \eqref{eq4.10} and using Lemma \ref{lem6.1}, we obtain the following corollary.

\bcor \label{cor4.x}
The process $D=D^r+D^g$ is a l\`agl\`ad process of finite variation and
\[
\xY=\xY_0+\wtf^r\bigcdot D^r+F^g \star D^g-(\xR-\xY)\bigcdot \Gamma-\xY_{-}G_{-}^{-1}\bigcdot \wt m+G^{-1}_{-}\,\bigcdot \wt\xN(\tau).
\]
\ecor

Assumption \ref{ass3.1}(vi) yields the existence $\FF$-predictable processes $\psi^{Y,Z}$ and $\nu$ such that
\begin{equation}  \label{eq4.11}
\xNX(\tau )=\psi^{\xY,\xZ}\bigcdot M, \quad m= \nu \bigcdot M.
\end{equation}
The next proposition is an immediate consequence of Corollary \ref{cor4.x} combined with \eqref{eq4.11}.
As before, we write $F_s^r(\cdot):=F_s^r(\cdot,\xZ_s,\xU_s)$ and $F_s^g(\cdot):= F_s^g(\cdot,\xZ_s,\xU_s)$
and we give an explicit representation for the $\FF$ BSDE associated with the $\GG$ BSDE \eqref{eq3.2}.

It is worth noting that Proposition \ref{pro4.1} extends several results from the existing literature where the method of reduction was studied
in a particular framework and under additional assumptions, such as the immersion hypothesis or the simplifying conditions  {\bf (A)} or  {\bf (C)}.

\bp  \label{pro4.1}
If the triplet $(\Ygg,\Zgg,\Ugg)$ is a solution to the $\GG$ BSDE \eqref{eq3.2}, then
the triplet $(\xY,\xZ,\xU)$ where $\xU = \Ugg$ satisfies on $\llb 0, \tau \rrb$
\begin{align*}
\xY_t&=\xX_{\tau}-\medint\int_{\,\rrb t,\tau\rrb} F^r_s(\xY_s)\,dD^r_s-\medint\int_{\,\llb t,\tau\llb} F_s^g (\xY_s)\,dD^g_{s+}-\medint\int_{\,\rrb t,\tau\rrb}z_s\,d\wt{M}_s\\
&\quad +\medint\int_{\,\rrb t,\tau\rrb}\big[\xR_s-\xY_s-(F^r_s(\xR_s)-F^r_s(\xY_s))\Delta D_s^r \big]\, d\Gamma_s
\end{align*}
where the process $z$ is given by $z_t:=G^{-1}_{t-}\big(\psi^{Y,Z}_t-Y_{t-}\nu_t\big)$.
\ep

In the above representation of the $\FF$ BSDE associated with the $\GG$ BSDE, we can clearly identify the reduced generators $F^r$ and $F^g$ and the form of the adjustment, which is integrated with respect to the hazard process $\Gamma= \wtG^{-1}\bigcdot A^o$ of a random time $\xtheta $.

\section{Construction of a solution to the $\GG$ BSDE}  \label{sect7}

In this section, we no longer postulate that a solution to the BSDE \eqref{eq3.2} exists, which means
that Assumption \ref{ass4.1} is relaxed. Our goal is to show that a solution to the $\FF$ BSDE
\eqref{eq5.6} can be expanded to obtain a solution to the BSDE \eqref{eq3.2} if equation \eqref{eq5.7}
has an $\FF$-optional solution $\xu$. We stress that equations \eqref{eq5.6} and \eqref{eq5.7} are
coupled, in the sense that they need to be solved jointly in order to construct a solution to the BSDE \eqref{eq3.2}.
Obviously, the issues of existence and uniqueness of a solution $(\xy,\xz,\xu)$ to equations \eqref{eq5.6} and \eqref{eq5.7}
need to be studied under additional assumptions on the generator and all other inputs to the BSDE \eqref{eq3.2}.
In the next result, we denote by $\ugg$ an arbitrary prescribed $\FF$-optional process and we do not use Lemma \ref{lem4.2}.

\bl \label{lem5.1}	
For a given process $\ugg \in \cOff$, let $(\xy,\xz )$ be an $\RR\times\RR^d$-valued, $\FF$-adapted solution to the BSDE, on $\llb 0, \tau \rrb $,
\begin{align} \label{eq5.1}
\xy_t& =\xX_{\tau}-\medint\int_{\,\rrb t,\tau\rrb} F_s^r(\xy_s,\xz_s,\widehat{\xu}_s)\,dD^r_s - \medint\int_{\,\llb t,\tau\llb} F_s^g(\xy_s,\xz_s,\widehat{\xu}_s)\,dD^g_{s+}  -\medint\int_{\,\rrb t,\tau\rrb}\xz_s\,d\wt{M}_s  \nonumber \\
&\quad +\medint\int_{\,\rrb t,\tau\rrb}\big[\xR_s-\xy_s-\big(F^r_s(\xR_s,\xz_s,\widehat{\xu}_s)
-F^r_s(\xy_s,\xz_s,\widehat{\xu}_s)\big)\Delta D_s^r \big]\dGammas
\end{align}
and let the $\GG$-adapted process $\ygg$ be given by
\begin{equation} \label{eq5.2}
\ygg:=\xy_0 + \I_{\rrb 0,\xtheta \llb} \bigcdot \xy^r + \I_{\rrb 0,\xtheta \llb} \star \xy^g +(\xR_\xtheta-\xy_{\xtheta-})\I_{\llb \xtheta ,\infty \llb}\I_\seq{\tau \geq \xtheta}.
\end{equation}
Then $(\ygg, \zgg ) := (\ygg,\xz^\xtheta )$ is a $\GG$-adapted solution to the BSDE, on $\llb 0, \tau\wedge\xtheta \rrb $,
\begin{align} \label{eq5.3}
\ygg_t &=X_{\tau\wedge \xtheta}-\medint\int_{\,\rrb t,\tau \wedge\xtheta\rrb }\wh F^r_s(\ygg_s,\zgg_s, \widehat{\xu}_s)\, dD^r_s \nonumber\\
&\quad -\medint\int_{\,\llb t,\tau \wedge\xtheta\llb }\wh F^g_s(\ygg_s,\zgg_s, \widehat{\xu}_s)\,dD^g_{s+}-\medint\int_{\,\rrb t,\tau \wedge\xtheta\rrb }\zgg_s\,d\wt{M}_s  \\
&\quad -\medint\int_{\,\rrb t,\tau \wedge\xtheta\rrb }\big[\xR_s-\xy_s-\big(F^r_s(\xR_s,\zgg_s,\ugg_s)-F^r_s(\xy_s,\zgg_s,\ugg_s)\big)\Delta D_s^r \big]d\Nogg_s .\nonumber
\end{align}
\el

\begin{proof}
Let us write ${C}^r:=\wtf^r \bigcdot D^r$ and $C^g:=F^g \star D^g_{+}$. From \eqref{eq5.1} and \eqref{eq5.2}, we can deduce that the
following equalities hold
\[
\ygg_{\tau \wedge\xtheta}= \xX_{\tau}\I_\seq{\tau<\xtheta}+R_{\xtheta}\I_\seq{\tau \geq \xtheta}
\]
and
\begin{align}  \label{eq5.4}
\ygg&= \xy_0+ \xz \I_{\rrb 0,\tau\rrb }\I_{\rrb 0,\xtheta \llb}\bigcdot\wt{M}^\xtheta -(\xR-\xy)\I_{\rrb 0,\tau\rrb}\I_{\rrb 0, \xtheta \llb}\bigcdot \Gamma \nonumber \\
& \quad +\I_{\rrb 0,\tau\rrb}\I_{\rrb 0,\xtheta \llb}\bigcdot{C}^r + \I_{\llb 0,\tau\llb}\I_{\llb 0,\xtheta \llb}\star {C}^g +(\xR-\xy_{-})\I_{\rrb 0,\tau\rrb }\bigcdot A.
\end{align}
Using again \eqref{eq5.1}, we obtain
\[
(\xy-\xy_-)\I_{\rrb 0,\tau \rrb}\bigcdot A 	
=\xz \I_{\rrb 0,\tau \rrb}\I_{\llb \xtheta \rrb}\bigcdot\wt{M}^\xtheta -(\xR-\xy)\I_{\rrb 0,\tau \rrb}\I_{\llb \xtheta \rrb}\bigcdot \Gamma
+\I_{\rrb 0,\tau \rrb}\I_{\llb \xtheta \rrb}\bigcdot{C}^r.
\]
Thus, by replacing $(R- \xy_-)$ by $(R-\xy)$ in the last term of \eqref{eq5.4} and using the equality
$\Nogg=A-\I_{\rrb 0,\xtheta\rrb} \bigcdot \Gamma $ (see \eqref{xeq2.1}), we see that $\ygg$ is equal to
\[
\xy_0 + \xz \I_{\rrb 0,\tau \wedge\xtheta\rrb}\bigcdot \wt{M}+(\xR-y)\I_{\rrb 0,\tau \wedge\xtheta\rrb }\bigcdot \Nogg+\I_{\rrb 0,\tau \wedge\xtheta\rrb} \bigcdot {C}^r + \I_{\llb 0,\tau \wedge\xtheta\llb}\star {C}^g_+.
\]
To establish \eqref{eq5.3}, it now remains to show that
\begin{equation}  \label{eq5.5}
\I_{\rrb 0,\tau \wedge\xtheta\rrb}\bigcdot C^r =\wh F^r(\ygg)\I_{\rrb 0,\tau \wedge\xtheta\rrb }\bigcdot D^r
-\big(F^r(R)-F^r(\xy)\big) \Delta D^r \I_{\rrb 0,\tau \wedge\xtheta\rrb }\bigcdot \Nogg
\end{equation}
where the variables $\xz$ and $\ugg$ are suppressed.

To this end, by using the fact that the equalities $F^r=\wh F^r$ and $\xy=\ygg$ hold on $\rrb 0,\xtheta\rrb$ and $\ygg \I_{\llb \xtheta \rrb}=R \I_{\llb \xtheta \rrb}$, we deduce from \eqref{eq4.4} that
\begin{align*}
I_1&:=\wtf^r(\xy,\widehat{\xu})\I_{\rrb 0,\tau\rrb}\I_{\rrb 0,\xtheta \llb}\bigcdot D^r=\wh F^r(\ygg)\I_{\rrb 0,\tau\wedge \xtheta\rrb}\bigcdot D^r  -F^r(R)\I_{\rrb 0,\tau\rrb}\I_{\llb \xtheta \rrb}\bigcdot D^r\\
& \quad +\Delta \Gamma \big(F^r(R)-F^r(\xy)\big)\I_{\rrb 0,\tau\rrb}\I_{\rrb 0, \xtheta \llb}\bigcdot D^r
\end{align*}
and
\[
I_2 :=\wtf^r (\xy)\I_{\rrb 0,\tau \rrb}\I_{\llb \xtheta \rrb}\bigcdot D^r=\big[ F^r(\xy)
+\Delta \Gamma \big(F^r(R)-F^r(\xy)\big)\big] \I_{\rrb 0,\tau \rrb}\I_{\llb \xtheta \rrb}\bigcdot D^r.
\]
Consequently, since  $\Nogg=A-\I_{\rrb 0,\xtheta\rrb}\bigcdot \Gamma$ is a process of finite variation stopped at $\xtheta $ and $\I_{\llb \xtheta \rrb}=\Delta A$, we conclude that
\begin{align*}
&\I_{\rrb 0,\tau \wedge\xtheta\rrb}\bigcdot {C}^r =I_1+I_2\\
& =\wh F^r(\ygg)\I_{\rrb 0,\tau \wedge\xtheta\rrb}\bigcdot D^r-\I_{\rrb 0,\tau\rrb}\big(F^r(R)-F^r(\xy)\big)
\big(\Delta A- \Delta \Gamma \I_{\rrb 0,\xtheta \rrb} \big) \bigcdot D^r\\
&=\wh F^r(\ygg)\I_{\rrb 0,\tau \wedge\xtheta\rrb}\bigcdot D^r-\I_{\rrb 0,\tau\rrb}\big(F^r(R)-F^r(\xy)\big)
\Delta D^r \bigcdot \big(A-\Delta \Gamma \I_{\rrb 0,\xtheta \rrb}\big)\\
&=\wh F^r(\ygg)\I_{\rrb 0,\tau\wedge\xtheta\rrb}\bigcdot D^r-\big(F^r(R)-F^r(\xy)\big)\Delta D^r\,\I_{\rrb 0,\tau \wedge\xtheta\rrb}\bigcdot \Nogg,
\end{align*}
which shows that equality \eqref{eq5.5} is valid.
\end{proof}

\brem
In the case where $R$ belongs to the class $\cPff$, one can modify the above proof by noticing that $Q:= A^o-A^p$ is a finite variation $\FF$-martingale
\[
R\wt G^{-1}\bigcdot (A^p-A^o)= -R\big(G_{-}^{-1} \bigcdot Q + \wt G^{-1}\bigcdot Q - G_{-}^{-1}\bigcdot Q\big) = -R G_{-}^{-1} \bigcdot \wt Q,
\]
which, when stopped at $\xtheta$, is a $\GG$-martingale. Then, in view of the predictable representation property of $M$,
this term will contribute to the $\GG$-martingale term $\widetilde M$ in \eqref{eq5.1}.
\erem

The following proposition is a consequence of Lemma \ref{lem5.1}. It shows that a solution to the $\GG$ BSDE can be constructed
by first solving the coupled equations \eqref{eq5.6}--\eqref{eq5.7}. Recall that we denote $F_s^r(\cdot):=F_s^r(\cdot,\xZ_s,\xU_s),\, F_s^g(\cdot):=  F_s^g(\cdot,\xZ_s,\xU_s),\,\wh F_s^r(\cdot ):=\wh F_s^r(\cdot,\Zgg_s,\Ugg_s)$ and $\wh F_g^r(\cdot ):=\wh F_s^g(\cdot,\Zgg_s,\Ugg_s).$

\bp \label{pro5.1}
Assume that $(\xy,\xz,\xu)$ is a solution to the BSDE on $\llb 0,\tau \rrb$
\begin{align}  \label{eq5.6}
\xy_t& =\xX_{\tau}-\medint\int_{\,\rrb t,\tau\rrb}F^r _s(\xy_s)\,dD^r_s -
\medint\int_{\,\llb t,\tau\llb} F^g _s(\xy_s)\,dD^g_{s+}-\medint\int_{\,\rrb t,\tau\rrb}\xz_s\,d\wt{M}_s\nonumber \\
& \quad +\medint\int_{\,\rrb t,\tau\rrb}\big[\xR_s-\xy_s- \big( F^r_s(\xR_s)-F^r_s(\xy_s) \big)\Delta D_s^r\big] \dGammas
\end{align}
where the $\FF$-optional process $\xu$ satisfies the following equality, for all $t \in \RR_+$,
\begin{equation}  \label{eq5.7}
\medint\int_{\,\rrb 0,t\rrb}\xu_s\,d\Nogg_s=\medint\int_{\,\rrb 0,t\rrb }\big[\xR_s-\xy_s-\big(F^r_s(\xR_s)-F^r_s(\xy_s)\big)\Delta D_s^r\big] \,d\Nogg_s .
\end{equation}
Then the triplet $(\ygg,\zgg ,\ugg ):= (\ygg, \xz^\xtheta , \xu )$ where the process $\ygg$ is given by
\[
\ygg:=\xy_0 + \I_{\llb 0,\xtheta \llb} \bigcdot \xy^r + \I_{\llb 0,\xtheta \llb} \star \xy^g +(\xR_\xtheta-\xy_{\xtheta-})\I_{\llb \xtheta ,\infty \llb}\I_\seq{\tau \geq \xtheta}
\]
is a solution to the BSDE \eqref{eq3.2} on $\llb 0, \xtheta \wedge \tau \rrb $, that is,
\begin{align*}
\ygg_t& =\Xgg_{\tau\wedge\xtheta}-\medint\int_{\,\rrb t,\tau\wedge\xtheta\rrb}\wh F^r_s(\ygg_s)\,dD_s^r - \medint\int_{\,\llb t,\tau\wedge\xtheta\llb} \wh F^g_s(\ygg_s)\,dD^g_{s+} \\
& \quad -\medint\int_{\,\rrb t,\tau\wedge\xtheta\rrb} \zgg_s\,d\wt{M}_s-\medint\int_{\,\rrb t,\tau\wedge\xtheta\rrb} \ugg_s \,d\Nogg_s.
\end{align*}
\ep

In the next step, we will examine the existence of a solution to the coupled equations \eqref{eq5.6}--\eqref{eq5.7}.
Specifically, we seek a triplet $(Y,Z,U)$ of processes that satisfy, for all $t \in [0, \tau]$,
\begin{align}  \label{eq39.0}
\xy_t&=\xX_{\tau}-\medint\int_{\,\rrb t,\tau\rrb}F_s^r(\xy_s)\,dD^r_s-\medint\int_{\,\llb t,\tau\llb} F_s^g(\xy_s)\,dD^g_{s+}-\medint\int_{\,\rrb t,\tau\rrb}\xz_s\,dM_s
\nonumber \\&\quad + \medint\int_{\,\rrb t,\tau\rrb}\xz_s\, G_s^{-1} \nu_s\, d[M,M]_s +\medint\int_{\,\rrb t,\tau\rrb}\big[\xR_s-\xy_s-\big(F^r_s(\xR_s)-F^r_s(\xy_s)\big)\Delta D_s^r \big]\dGammas
\end{align}
and, for all $t \in \RR_+$ (notice that \eqref{eq39.1} is manifestly stronger than \eqref{eq5.7})
\begin{equation} \label{eq39.1}
\xu_t=\xR_t-\xy_t- \big(F^r_t(\xR_t)-F^r_t(\xy_t)\big)\Delta D_t^r.
\end{equation}

To examine the existence of a solution to the coupled equations  \eqref{eq39.0}--\eqref{eq39.1}, we introduce the transformed equations \eqref{eq40.0}--\eqref{eq40.1}. Our goal is to remove the quadratic variation term $G^{-1} \bigcdot [m,M] = G^{-1} \nu \bigcdot [M,M]$ from \eqref{eq39.0} and place $\nu$ inside the generators $F^r$ and $F^g$, which are assumed to be bounded.
In that way, we avoid the need to check the appropriate growth conditions when applying the existing results on well-posedness of BSDEs.

We define the linear transformation $\bar \xy := G\xy$, $\bar \xz:=G_{-}\xz+G^{-1}\bar \xy\nu$ and $\bar \xu:=\xu$. Hence we obtain the transformed generators
\begin{align*}
\bar F^{r}_s(y,z,u) := \wt G_s F^r_s(G_s^{-1}y,G_{s-}^{-1}(z- G^{-1}_sy\nu_s),u)
\end{align*}
and
\begin{align*}
\bar F^{g}_s(y,z,u) := G_s F^g_s(G_s^{-1}y,G_{s-}^{-1}(z- G^{-1}_sy\nu_s),u)
\end{align*}
and we denote $\bar F_s^r(\cdot ):=\bar F_s^r(\cdot ,\bar \xz_s,\bar \xu_s)$ and $\bar F_s^g(\cdot ):=\bar F_s^g(\cdot ,\bar \xz_s,\bar \xu_s)$. Observe that if a solution $(\xy, \xz,\xu) \in \cOff \times \cPffd \times \cOff$ to \eqref{eq39.0}--\eqref{eq39.1} exists, then $(\bar \xy, \bar \xz, \bar \xu) \in \cOff \times \cPffd \times \cOff$ satisfies the following coupled equations, for all $t \in [0, \tau]$,
\begin{align} \label{eq40.0}
\bar \xy_t& =G_\tau\xX_{\tau}-\medint\int_{\,\rrb t,\tau\rrb} \bar F_s^r(\bar \xy_s)\,dD^r_s - \medint\int_{\,\llb t,\tau\llb} \bar F_s^g(\bar \xy_s)\,dD^g_{s+}\nonumber \\
& \quad  -\medint\int_{\,\rrb t,\tau\rrb} \bar \xz_s\,dM_s +\medint\int_{\,\rrb t,\tau\rrb}\big[\wt G_s\xR_s-\Delta \bar F^r_s\Delta D_s^r \big]\,d\Gamma_s
\end{align}
 and, for all $t\in \RR_+$
\begin{equation} \label{eq40.1}
\bar \xu_t = \xR_t-\bar \xy_sG^{-1}_t- \wt G^{-1}_t\Delta \bar F^r_t  \Delta D_t^r
\end{equation}
where we denote $\Delta \bar F^r_s := \bar F^r_s(G_s\xR_s)-\bar F^r_s(\bar \xy_s)$. In the reverse, a solution $(\xy, \xz, \xu)$ to the coupled equations \eqref{eq39.0}--\eqref{eq39.1} can be obtained from a solution $(\bar \xy, \bar \xz ,\bar \xu)$ to the coupled equations \eqref{eq40.0}--\eqref{eq40.1} by setting $\xy := G^{-1}\bar \xy$, $Z := G^{-1}_-(\bar \xz - G^{-1} \bar \xy \nu)$ and $\xu := \bar \xu$.

Observe that BSDEs \eqref{eq39.0} and \eqref{eq40.0} have the l\`agl\`ad driver $D=(D^r,D^g)$, which may share common jumps with the martingale $M$. To the best of our knowledge, there is a gap in the literature on BSDEs of this form and thus we develop in Section \ref{sect12} a jump-adapted methodology to solve such BSDEs under specific assumptions.

\brem
If a random time $\xtheta$ is an $\FF$-{\it pseudo-stopping time} (see Nikeghbali and Yor \cite{NY2005} and Aksamit and Li \cite{AL}), then $m=1$ and thus $\nu= 0$ in \eqref{eq4.11} so that we can deal directly with \eqref{eq39.0}. Also, under Assumption {\bf (C)}, it is possible to eliminate the term $G^{-1} \bigcdot [m,M]$ using the Girsanov transform in $\GG$ if $\xtheta$ is an {\it invariance time}, as introduced and studied in Cr\'epey and Song \cite{CS2017}. We point out that in \cite{CS2017}, the authors have worked under a predictable setup, rather than an optional one, as we do in the present work. Furthermore, the $\GG$ BSDE in \cite{CS2015} was stopped strictly before $\xtheta$ (specifically, at $\xtheta-$), rather than at $\xtheta$.
Therefore, the Girsanov theorem was in fact applied  in \cite{CS2015} to eliminate the term $G^{-1} \bigcdot [n,M]$, which is known to coincide with $G^{-1} \bigcdot [m,M]$ under Assumption {\bf (C)}. In principle, one could attempt to mimic the approach developed in \cite{CS2017} by studying a new class of random times for which the drift term $G^{-1} \bigcdot [m,M]$ can be removed or, in other words, to introduce a family of random times for which there exists an equivalent probability measure under which they are $\FF$-pseudo-stopping times. However, this issue is beyond the scope of the current work and here we assume that either the generator is bounded or $\xtheta$ is an $\FF$-pseudo-stopping time.
\erem

\section{Solution to the $\GG$ BSDE in the Brownian case} \label{sect8}

To illustrate our approach from Sections \ref{sect6} and \ref{sect7}, we first use results from Essaky et al. \cite{EHO2015} to show that in the case of the Brownian filtration, if $F^g = 0$ and $U$ does not appear in the right-hand side of \eqref{eq39.1} or \eqref{eq40.1}, then a unique solution $(\bar \xy, \bar \xz,\bar \xu)$ to \eqref{eq40.0}--\eqref{eq40.1} exists and thus a unique solution to $(\xy, \xz,\xu)$ to \eqref{eq39.0}--\eqref{eq39.1} exists as well.
Specifically, we now take $M=W$ to be a $d$-dimensional Wiener process in its natural filtration $\FF$ and we consider below the coupled equations \eqref{eq39.0}--\eqref{eq39.1} with $F^g= 0$. Our goal is to demonstrate that, under some natural assumptions, these equations have
a unique solution $(\xy,\xz,\xu) \in \cOff \times \cPffd \times \cOff$.

To simplify the notation, we write $D$ and $F$ instead of $(D^r,0)$ and $(F^r,0)$, respectively, and we consider the situation where $D=(D^1,D^2)=(\langle W \rangle,\Gamma )$ and $F=(F^1,F^2)=(F^1({y},{z},{u}),F^2({y}))$ for $(y,z,u) \in \RR\times \RR^d\times \RR$. Then the BSDE \eqref{eq40.0} becomes
\begin{align}
\bar \xy_t &=G_\tau\xX_{\tau}-\medint\int_{\,\rrb t,\tau\rrb} \bar F_s^1(\bar \xy_s)\,d\langle W\rangle_s -\medint\int_{\,\rrb t,\tau\rrb} \bar F_s^2(\bar \xy_s)\,d\Gamma_s -\medint\int_{\,\rrb t,\tau\rrb} \bar \xz_s\,d W_s \nonumber \\	
&\quad +\medint\int_{\,\rrb t,\tau\rrb}\big[\wt G_{s}\xR_s-\big( \bar F^2_s(G_s\xR_s)-\bar F^2_s(\bar \xy_s) \big)\Delta \Gamma_s \big]\,d\Gamma_s \label{eq43}
\end{align}
where $\bar F_s^1(\bar \xy_s):=\bar F_s^1(\bar \xy_s,\bar \xz_s,\bar \xu_s)$ and equation \eqref{eq40.1} has an explicit solution given by
\begin{equation} \label{eq43.0}
\bar \xu = (\xR-\bar \xy G^{-1}) - (\bar F^2(G\xR)-\bar F^2(\bar \xy) )\wt G^{-1} \Delta \Gamma .
\end{equation}
We point out that, in the above, our choice of $D^2$ was somewhat arbitrary and we decided to set $D^2=\Gamma$ for simplicity of presentation.

\bp \label{pro9.1}
Assume that:\\
(i) for every $t\in \RR_+$, the map $F^1_t:\RR\times \RR^{d}\times\RR\to\RR$ is bounded and Lipschitz continuous;\\
(ii) for every $t\in \RR_+$, the map $F^2_t:\RR\to\RR$ is bounded, Lipschitz continuous and decreasing;\\
(iii) the dual $\FF$-optional projection $A^o$ has a finite number of discontinuities.\\
Then the BSDE \eqref{eq43} has a solution $(\bar \xy,\bar \xz)$ and a solution $(\xy,\xz)$ to \eqref{eq39.0} is obtained by
setting $\xy:=G^{-1}\bar \xy$ and $\xz:= G^{-1}(\bar \xz - G^{-1}\bar \xy \nu )$.
\ep

\proof
To establish the existence of a solution $(\bar\xy,\bar \xz)$ to \eqref{eq43}, we will apply Theorem 2.1 in Essaky et al. \cite{EHO2015} to the  data $(P, R,\xtheta,F)$. We observe that the BSDE \eqref{eq43} is a special case of equation (2.1) in \cite{EHO2015}, which has the following form
\begin{align} \label{eq43.1}
\bar Y_t&=G_\tau\xX_\tau+\medint\int_{\,\rrb t,\tau\rrb}f(s,\bar Y_s,\bar Z_s)\,ds+\medint\int_{\,\rrb t,\tau\rrb}g(s,\bar Y_s)\,d\bar{A}_s \nonumber \\
&\quad +\sum_{t<s \leq \tau}h(s,\bar Y_{s-},\bar Y_s)-\medint\int_{\,\rrb t,\tau\rrb}\bar Z_s\,dW_s .
\end{align}

Indeed, \eqref{eq43} can be recovered from \eqref{eq43.1} if we set $\bar{A}:=\Gamma^c$ where $\Gamma^c$ is the continuous part of $\Gamma$ and
define the mappings $f,g$ and $h$ as follows
\begin{align*}
& f(s,\bar \xy_s,\bar \xz_s):=-\bar F^1_s(\bar Y_s,\bar Z_s,\bar \xu_s), \quad g(s,\bar \xy_s):=\wt G_s\xR_s-\bar F^2_s(\bar \xy_s),\\
& h(s,\bar \xy_{s-},\bar \xy_s):=\big(\wt G_s\xR_s-F^2_s(\bar \xy_s)\big)\Delta \Gamma_s
-\big(\bar F^2_s(G_s \xR_s)-\bar F^2_s(\bar \xy_s)\big)(\Delta \Gamma_s)^2
\end{align*}
where $\bar \xu$ is given by \eqref{eq43.0} and the mapping $h$ is in fact independent of the variable $\bar Y_{-}$.

To apply Theorem 2.1 from Essaky et al. \cite{EHO2015}, it suffices to check that Assumptions (A.1)--(A.4) on page 2151 of \cite{EHO2015} are satisfied.
To check Assumption (A.1), we note that the mapping $F^1$ is assumed to be continuous and, since it is also bounded, there exists a constant $C>0$ such that $|F^1|\leq C \leq C(1+|z|)$, which shows that the linear growth condition is satisfied. Next, from the fact that $|\Delta \Gamma| \leq 1$, we deduce that $|g|$ is bounded and thus Assumption (A.2) holds as well. Finally, $h$ is bounded and continuous and, since the process $\Gamma$ is assumed to have a finite number of jumps, it is enough to check condition (c) in Assumption (A.3).

To this purpose, we observe that $0\leq \Delta \Gamma \leq 1$ and thus the mapping
\[
y \mapsto y+ h(s,y)=y-\Delta \Gamma_s \bar F^2_s(y)\big(1-\Delta \Gamma_s\big)+\Delta\Gamma_s\big(\wt G_s \xR_s-\Delta \Gamma_s \bar F^2_s(G_s\xR_s)\big)
\]
is nondecreasing and continuous. Finally, the Mokobodski condition postulated in (A.4) is trivially satisfied as we
deal here with the BSDE with no reflecting boundaries. Thus, by applying Theorem 2.1 in \cite{EHO2015} with $T$ replaced
by $\tau \in\cT$, we obtain the existence of a maximal solution $(\bar \xy, \bar \xz)$ to \eqref{eq43}.
\endproof

\bex \label{ex9.1}
Let us consider a special case where $\xtheta$ is an $\FF$-pseudo-stopping time (see Nikeghbali and Yor \cite{NY2005}) 
and thus the process $\nu$ in \eqref{eq4.11} vanishes. We can suppose that the mapping $z\mapsto F^1_s(y,z,u)$ has a linear (quadratic) growth and thus, since in the case of a Brownian filtration the equality $\wt G = G_-$ holds, we obtain
\[
|\bar F^1_s(y,z)|=\wt G_s\big|F^1_s(G_s^{-1}y,G_{s-}^{-1}(z- G^{-1}_sy\nu_s),u)\big| \leq G_{s-}(1+G^{-1}_{s-}|z|)\leq (1 + |z|),
\]
which shows that the boundedness of $F^1$ postulated in Proposition \ref{pro9.1} can be relaxed.
\eex

\section{Reduction of a solution to the $\GG$ RBSDE}     \label{sect9}

Our goal in this section is to study the properties of solutions to the $\GG$ RBSDE with a random time horizon $\xtheta $.

\bd \label{def3.4}
A quadruplet $(\Ygg,\Zgg,\Ugg, \Lgg)$ is a {\it solution} on the interval $\llb 0, \tgg \wedge \xtheta\rrb $ to the $\GG$ RBSDE
\begin{align} \label{eq3.5}
\Ygg_t&=\Xgg_{\tgg\wedge\xtheta}-\medint\int_{\,\rrb t,\tgg\wedge\xtheta\rrb}\wh F_s^r(\Ygg_s,\Zgg_s,\Ugg_s)\,d\wh D^r_s-\medint\int_{\,\llb t,\tgg\wedge \xtheta\llb}
\wh F_s^g(\Ygg_s,\Zgg_s,\Ugg_s)\,d\wh D^g_{s+} \nonumber \\
&\quad -\medint\int_{\,\rrb t,\tgg\wedge\xtheta\rrb}\Zgg_s\,d\wt{M}_s-\medint\int_{\,\rrb t,\tgg\wedge \xtheta\rrb}\Ugg_s\,d\Nogg_s-(\Lgg_{\tgg\wedge\xtheta}-\Lgg_{t})
\end{align}
if $\Ygg\in\cOgg$ is a l\`agl\`ad process such that $\Ygg \geq \Xgg $, the processes $\Zgg\in\cPggd$ and $\Ugg\in\cOff$
are such that the stochastic integrals in the right-hand side of \eqref{eq3.2} are well defined, $\Lgg = \Lgg^\xtheta$ is a
l\`agl\`ad, increasing, and $\GG$-strongly predictable process with $\Lgg_0= 0$ and has the decomposition $\Lgg=\Lgg^r+\Lgg^g$ where the processes $\Lgg^r$ and $\Lgg^g$ obey the Skorokhod conditions
\[
\big(\I_{\{\Ygg_{-}\neq \Xgg_{-}\}}\bigcdot \Lgg^r\big)_{\tgg}^\xtheta =\big(\I_{\{\Ygg\neq \Xgg\}}\star \Lgg^g_{+}\big)_{\tgg}^\xtheta=0
\]
and equality \eqref{eq3.5} is assumed to hold on $\llb 0, \tgg \wedge \xtheta\rrb$.
\ed

As in Section \ref{sect6}, in order to show that \eqref{eq3.5} has a solution, we will first reduce a solution to
the $\GG$ RBSDE \eqref{eq3.5} to a solution of an associated RBSDE in filtration $\FF$. Subsequently, we show how a solution
to the reduced $\FF$ RBSDE, which can be shown to exist under suitable assumptions, can be employed to construct a solution to the $\GG$ RBSDE \eqref{eq3.5}.
Again, we first work under the following temporary postulate, which will be relaxed in Section \ref{sect10}.

\bhyp \label{ass6.1}
A solution $(\Ygg,\Zgg,\Ugg,\Lgg)$ to the $\GG$ RBSDE \eqref{eq3.5} exists.
\ehyp

Following the approach developed for the non-reflected case, we decompose $\Ygg$ into the pre-default and post-default components
\begin{align*}
\Ygg_t\I_{\{\xtheta >t\}}+\Ygg_t\I_{\{\xtheta \leq t\}}
&=\Ygg_t\I_{\{\xtheta > t\}}+\Xgg_{\xtheta \wedge \tau}\I_{\{\xtheta \leq t\}}=Y_t\I_{\{\xtheta >t\}}+\Xgg_{\xtheta}\I_{\{\xtheta \leq t\}}\\
&=G^{-1}_t\EE\big(\Ygg_t\I_{\{\xtheta >t\}}\,\big|\,\cF_t\big)\I_{\{\xtheta >t\}}+R_{\xtheta}\I_{\{\xtheta \leq t\}}.
\end{align*}
To compute the component $G^{-1}_t\EE\big(\Ygg_t\I_{\{\xtheta >t\}}\,|\,\cF_t\big)$ we proceed similarly to the non-reflected case. The new feature here is the use of Lemma \ref{lem6.1} in order to obtain a reduction of the $\GG$-strongly predictable, increasing process $\Lgg$.
Computations in this section are similar to those in Section \ref{sect6}, except for the presence of the reflection process
and thus in the following we will focus on new elements. To proceed, similarly to \eqref{eq4.8}, we set
\begin{align} \label{K}
K_t(\tau)&:= \EE\big(\xX_{\tau}G_{\tau} +(\wt G \wtf^r\bigcdot D)_\tau +(GF^g\star D^g_+)_{\tau}\,|\,\cF_t\big) \nonumber \\
& \quad  +  \EE \big((\xR+\xL)\bigcdot A^o)_\tau +\xL_\tau G_\tau\,|\,\cF_t\big)
\end{align}
where $\wtf^r$ is given by \eqref{eq4.4}. From Assumption \ref{ass3.1}, we deduce the existence of $\FF$-predictable processes
$\psi^{\xY,\xZ,\xL}$ and $\nu$ such that $K(\tau)=\psi^{\xY,\xZ,\xL}\bigcdot M $ and $m = \nu \bigcdot M.$

\bp \label{pro6.1}
The process $\xY = \,^o(\Ygg\I_{\llb 0, \xtheta \llb})G^{-1}$ satisfies the $\FF$ RBSDE, on $\llb 0,\tau \rrb $,
\begin{align*}  
\xY_t & =\xX_\tau-\medint\int_{\,\rrb t,\tau\rrb}F_s^r(\xY_s)\,dD^r_s-\medint\int_{\,\llb t,\tau \llb}F_s^g(\xY_s)\,dD^g_{s+}-\medint\int_{\,\rrb t,\tau \rrb } z_s\,d\wt M_s\nonumber \\ &\quad + \medint\int_{\,\rrb t,\tau\rrb }\big[\xR_s-\xY_s -(F^r_s(\xR_s)-F^r_s(\xY_s))\Delta D_s^r \big]\dGammas
-(\xL_\tau-\xL_t)
\end{align*}
where $\xY \geq \xX,\, \xU = \Ugg$, the process $\xZ$ is the $\FF$-reduction of the process $\Zgg$ given by \eqref{eq4.1}, $z_t:=G^{-1}_{t-}\big(\psi^{\xY,\xZ,\xL}_t-\xY_{t-}\nu_t\big)$ is an $\FF$-predictable process and $\xL$ is an $\FF$-strongly predictable, increasing process such that $\xL = \Lgg$ on $\llb 0, \xtheta\rrb$ and the Skorokhod conditions are satisfied, that is, $(\I_{\{\xY_-\neq \xX_-\}}\bigcdot \xL^r)_\tau=(\I_{\{\xY\neq \xX\}}\star \xL^g_+)_{\tau}=0$.
\ep

\begin{proof}
In view of Lemmas \ref{lem4.3},  \ref{lem6.1} and \ref{lem4.2} the generators $\wh F^r$, $\wh F^g$ and the increasing process $\wh L$ can be reduced to the filtration $\FF$ to obtain, on the event $\{\tau \geq t\}$,
\begin{align*}
\EE&(\Ygg_t\I_{\{\xtheta>t\}}\,|\,\cF_t)
=\EE\big(P_{\tau}\I_{\{\xtheta > \tau\}}+R_{\xtheta}\I_{\{t< \xtheta \leq \tau\}}+(\wt G \wtf^r\bigcdot D^r)_{\tau} - (\wt G \wtf^r\bigcdot D^r)_t
\\ & \quad + (Gf^g\bigcdot D^g_{+})_{\tau} - (Gf^g\bigcdot D^g_+)_t
 + (\xL_{\tau}-\xL_t)\I_{\{\xtheta > \tau\}} + (\xL_{\xtheta}-\xL_t)\I_{\{t< \xtheta \leq \tau\}} \,|\,\cF_t\big)\\
&  =\EE\big(\xX_{\tau}G_\tau + (\wt G \ddot  F^r\bigcdot D^r)_{\tau} - (\wt G \wtf^r\bigcdot D^r)_t + (Gf^g\bigcdot D^g_{+})_{\tau} - (Gf^g\bigcdot D^g_+)_t
\\& \quad + ((\xR+\xL)\bigcdot A^o)_\tau-((R+L)\bigcdot A^o)_t+\xL_\tau G_\tau- \xL_tG_t\,|\,\cF_t\big)
\end{align*}
where the mapping $\wtf^r$ is given in \eqref{eq4.4}.

Next, an application of the l\`agl\`ad product rule to $LG$ and the equalities $\wtG=G_- + \Delta m$ and $\xL=\xL_- + \Delta \xL^r$ yield
\begin{align*}
LG&=\xL_{-}\bigcdot m-\xL_{-} \bigcdot A^o+G_{-}\bigcdot L^r+G\star L^g_{+}+\Delta G \bigcdot \Delta \xL^r\\
  &=\xL_{-}\bigcdot m-\xL\bigcdot A^o+\wtG\bigcdot \xL^r+G\bigcdot  \xL^g_{+}.
\end{align*}
By combining these computations, we conclude that
\[
\xY G =K(\tau)- \wt G\wtf^r\bigcdot D^r - G F^g\star D^g_{+} - \xR\bigcdot A^o+\wtG\bigcdot \xL^r+ G\star \xL^g_{+}
\]
where $\xY:=G^{-1} \,^o(\Ygg\I_{\llb 0, \xtheta \llb})$ and $K(\tau)$ is given by \eqref{K}.

The backward dynamics of $\xY$ can now be computed from
Corollary \ref{cor4.x}
\begin{align*}
\xY_t & =\xX_\tau - \medint\int_{\,\rrb t,\tau\rrb }\wtf_s(\xY_s)\,dD^r_s - \medint\int_{\,\llb t,\tau \llb }F^g_s(\xY_s)\,dD^g_{s+}-\medint\int_{\,\rrb t,\tau\rrb } z_s\,dM_s  \\
	& \quad +\medint\int_{\,\rrb t,\tau\rrb }(\xR_s-\xY_s)\dGammas +\medint\int_{\,\rrb t,\tau\rrb }\wtG^{-1}_s z_s\,d[M,n]_s-(\xL_\tau-\xL_t)
\end{align*}
and thus, after rearranging and using \eqref{eq4.4}, we obtain the asserted BSDE.

It remains to check that the appropriate Skorokhod conditions are met by the process $\xL$. Recall that the Skorokhod conditions satisfied by $\Lgg^{r}$ and $\Lgg^{g}$ are
\begin{equation} \label{eq6.2}
\big(\I_{\{\Ygg_-\neq X_-\}}\bigcdot \Lgg^r\big)_{\tau\wedge \xtheta}=\big(\I_{\{\Ygg \neq X\}}\star \Lgg^g_+\big)_{\tau\wedge \xtheta} = 0.
\end{equation}
By integrating the first equality in \eqref{eq6.2} with respect to $G^{-1}_-$, we obtain
\[
\big(G^{-1}_-\I_{\llb 0,\xtheta \rrb}\I_{\{\Ygg_-\neq X_-\}}\bigcdot \Lgg^r\big)_\tau=0.
\]
The equality $\Ygg\I_{\llb 0,\xtheta \llb}=\xY\I_{\llb 0,\xtheta \llb}$ implies that $\Ygg_-\I_{\llb 0,\xtheta \rrb}=\xY_-\I_{\llb 0,\xtheta \rrb}$ and $(\Lgg^r)^\xtheta=(\xL^{r})^{\xtheta}$. Consequently, since $X_-=\xX_-\I_{\llb 0, \xtheta \rrb}+R_{\xtheta}\I_{\rrb \xtheta,\infty\rrb}$,
we get
\[
\big(G^{-1}_-\I_{\llb 0,\xtheta \rrb}\I_{\{\Ygg_-\neq X_-\}}\bigcdot \Lgg^r\big)_\tau =
\big(G^{-1}_-\I_{\llb 0,\xtheta \rrb}\I_{\{\xY_-\neq \xX_-\}}\bigcdot \xL^r\big)_\tau=0 .
\]
Then, by taking the expectation and using the property of the dual $\FF$-predictable projection, we obtain
\[
\EE( (G^{-1}_-\I_{\llb 0,\xtheta \rrb}\I_{\{\xY_-\neq \xX_-\}}\bigcdot \xL^r)_\tau) =\EE((\I_{\{\xY_-\neq \xX_-\}}\bigcdot \xL^r)_\tau),
\]
which implies that $\xL^r$ obeys the first Skorokhod condition, that is, $(\I_{\{\xY_-\neq \xX_-\}}\bigcdot \xL^r)_\tau=0$.
Similarly, to check the second Skorokhod condition, we integrate the second equality in \eqref{eq6.2} with respect to $G^{-1}$ and use the equality $\Lgg^g_+\I_{\llb 0 ,\xtheta \llb}=\xL^g_+\I_{\llb 0 ,\xtheta \llb}$ to obtain
\[
(G^{-1}\I_{\llb 0,\xtheta \llb}\I_{\{\Ygg\neq X\}}\star \Lgg^g_+)_{\tau}=(G^{-1}\I_{\llb 0,\xtheta \llb}\I_{\{\xY \neq P \}}\star \xL^{g}_+)_{\tau}.
\]
By taking the expectation and using the property of the dual $\FF$-optional projection, we obtain the equality $\EE((\I_{\{\xY \neq \xX \}}\star \xL^{g}_+)_{\tau})=0$, which in turn implies that $(\I_{\{\xY \neq \xX \}}\star \xL^{g}_+)_{\tau}=0$. 
\end{proof}

\section{Construction of a solution to the $\GG$ RBSDE} \label{sect10}

In this section, we relax the postulate that a solution $(\Ygg, \Zgg, \Ugg, \Lgg)$ exists.
In the next result, we again denote by $\wh {\xu }$ an arbitrary prescribed $\FF$-adapted process and we do not use Lemma \ref{lem4.2}.

\bl \label{lem7.1}	
Let a process $\widehat{\xu} \in \cOff $ be given and let $(\xy,\xz,\xl)$ be an $\RR\times\RR^d\times \RR$-valued, $\FF$-adapted solution to the $\FF$ RBSDE, on $\llb 0, \tau \rrb $,
\begin{align*}
\xy_t &=\xX_{\tau}-\medint\int_{\,\rrb t,\tau\rrb} F_s^r(\xy_s,\xz_s,\widehat \xu_s)\,dD_s^r- \medint\int_{\,\llb t,\tau\llb} F_s^g (\xy_s,\xz_s,\ugg_s)\,dD_{s+}^g\\
&\quad +\medint\int_{\,\rrb t,\tau\rrb}\big[\xR_s-\xy_s-\big(F_s^r(\xR_s,\xz_s,\widehat{\xu}_s)-F^r_s(\xy_s,\xz_s,\widehat{\xu}_s)\big)\Delta D_s^r \big]\dGammas\\
&\quad -\medint\int_{\,\rrb t,\tau\rrb}\xz_s\,d\wt{M}_s - (\xl_\tau- \xl_t)
\end{align*}
where $\xy \geq \xX$ and $\xl$ is an $\FF$-strongly predictable, increasing process such that the Skorokhod conditions $(\I_{\{\xy_-\neq \xX_-\}}\bigcdot \xl^r)_\tau=(\I_{\{\xy\neq \xX\}}\star  \xl^g_+)_{\tau}=0$ hold and $\xL_0 = 0$.
Then the triplet $(\ygg,\zgg ,\lgg ):= (\ygg,\xz^\xtheta ,\xl^\xtheta)$ where $\ygg$ is given by
\[
\ygg:= \xy_0 + \I_{\rrb 0,\xtheta \llb}\bigcdot \xy^r + \I_{\rrb 0,\xtheta \llb}\star \xy^g_{+}
+(\xR_\xtheta-\xy_{\xtheta-})\I_{\llb \xtheta ,\infty \llb}\I_\seq{\tau \geq \xtheta}
\]
is a solution to the $\GG$ RBSDE, on $\llb 0, \tau\wedge\xtheta \rrb $,
\begin{align*}
\ygg_t&=X_{\tau\wedge \xtheta} - \medint\int_{\,\rrb t,\tau \wedge\xtheta\rrb }\wh F^r_s(\ygg_s,\zgg_s,\widehat{\xu}_s)\, dD^r_s
- \medint\int_{\,\llb t,\tau \wedge\xtheta\llb }\wh F^g_s(\ygg_s,\zgg_s,\ugg_s)\, dD^g_{s+}\\
&\quad - \medint\int_{\,\rrb t,\tau \wedge\xtheta\rrb }\big[\xR_s-\xy_s-\big(F_s^r(\xR_s,\zgg_s,\ugg_s)-F^r_s(\xy_s,\zgg_s,\ugg_s)\big)\Delta D_s \big]\,d\Nogg_s \\
&\quad - \medint\int_{\,\rrb t,\tau \wedge\xtheta\rrb }\zgg_s\,d\wt{M}_s -(\lgg_{\tau\wedge \xtheta}- \lgg_t) \nonumber
\end{align*}
where $\ygg \geq \xX$ and $\lgg = \xl^\xtheta$ is a $\GG$-strongly predictable, increasing process such that the Skorokhod conditions  $(\I_{\{\ygg_-\neq \Xgg_-\}}\bigcdot \lgg^r)_{\xtheta \wedge \tau}=(\I_{\{\ygg \neq \Xgg \}}\star \lgg^g_+)_{\xtheta \wedge \tau}=0$ are valid and $\lgg_0 = 0$.
\el

\begin{proof}
It suffices to set  ${C}^r:=\wtG^{-1}\wtf^r\bigcdot D^r + \xl^r$ and ${C}^g:= F^g\star D_{+}^g + \xl^g$
in the proof of Lemma \ref{lem5.1}. The required Skorokhod conditions are also met since
\begin{align*}
\big(\I_{\llb 0,\xtheta \rrb}\I_{\{\ygg_-\neq X^{\xtheta}_-\}}\bigcdot \xl^{r}\big)_\tau   =
\big(\I_{\llb 0,\xtheta \rrb}\I_{\{\xy_-\neq \xX_-\}}\bigcdot \xl^r\big)_\tau=0
\end{align*}
and
\begin{align*}
\big(\I_{\llb 0,\xtheta \llb}\I_{\{\ygg\neq X^{\xtheta}\}}\star \xl^{g}_+\big)_\tau    =
\big(\I_{\llb 0,\xtheta \llb}\I_{\{ \xy\neq \xX\}}\star \xl^g_+\big)_\tau=0
\end{align*}
where we have used the following equalities: $\ygg\I_{\llb 0,\xtheta \llb}= \xy\I_{\llb 0,\xtheta \llb},\,\ygg_-\I_{\llb 0,\xtheta \rrb}=\xy_-\I_{\llb 0,\xtheta \rrb}$ and $\xX_-=\xX_-\I_{\llb 0, \xtheta \rrb}+R_{\xtheta}\I_{\rrb \xtheta,\infty\rrb}$.
\end{proof}

\bp \label{pro7.1}
Let $(\xy,\xz,\xu, \xl)$ be a solution to the $\FF$ RBSDE, on $\llb 0,\tau \rrb$,
\begin{align}
\xy_t&=\xX_\tau -\medint\int_{\,\rrb t,\tau\rrb }F^r_s(\xy_s)\,dD^r_s -\medint\int_{\,\llb t,\tau\llb }F_s^g(\xy_s)\,dD^g_{s+}-\medint\int_{\,\rrb t,\tau\rrb } \xz_s\,d\wt M_s   \label{eq53.0} \\
&\quad + \medint\int_{\,\rrb t,\tau\rrb }\big[\xR_s-\xy_s -\big(F^r_s(\xR_s) - F^r_s(\xy_s)\big)\Delta D_s^r\big]\dGammas
-(\xl_{\tau } - \xl_t) \nonumber
\end{align}
where $\xy \geq \xX,\, \xl$ is an $\FF$-strongly predictable, increasing process  with $L_0=0$ and such that the Skorokhod conditions $(\I_{\{\xy_-\neq \xX_-\}}\bigcdot \xl^r)_\tau=(\I_{\{ \xy\neq \xX \}}\star \xl^g_+)_{\tau}=0$ are obeyed, and the $\FF$-optional process $\xu$ satisfies, for all $t \in \RR_+$,
\begin{equation}  \label{eq53.1}
\medint\int_{\,]0,t]}\xu_s\,d\Nogg_s=\medint\int_{\,]0,t]}\big[\xR_s-\xy_s-\big(F^r_s(\xR_s)- F^r_s(\xy_s)\big)\Delta D_s^r\big]\,d\Nogg_s .
\end{equation}
Then $(\ygg,\zgg ,\ugg, \lgg ):= (\ygg,\xz^\xtheta ,\xu, \xl^\xtheta )$ where the process $\ygg$ is given by
\[
\ygg:=\xy_0 + \I_{\llb 0,\xtheta \llb}\bigcdot \xy^r + \I_{\llb 0,\xtheta \llb}\star \xy^g_{+} +(\xR_\xtheta-\xy_{\xtheta-})\I_{\llb \xtheta ,\infty \llb}\I_\seq{\tau \geq \xtheta}
\]
is a solution to the $\GG$ RBSDE \eqref{eq3.5} on $\llb 0,\tau\wedge\xtheta \rrb$, that is,
\begin{align} \label{prop3.2eq1}
\ygg_t&=X_{\tau\wedge \xtheta}-\medint\int_{\,\rrb t,\tau \wedge\xtheta\rrb }\wh F^r_s(\ygg_s)\, dD^r_s-\medint\int_{\,\llb t,\tau \wedge\xtheta\llb }\wh F^g_s(\ygg_s)\, dD^g_{s+}\\
&\quad -\medint\int_{\,\rrb t,\tau \wedge\xtheta\rrb }\zgg_s\,d\wt{M}_s-\medint\int_{\,\rrb t,\tau \wedge\xtheta\rrb }\ugg_s \,d\Nogg_s  -( \lgg_{\xtheta\wedge \tau}- \lgg_t ) \nonumber
\end{align}
where $\ygg \geq X$ and $\lgg = \xl^\xtheta$ is a $\GG$-strongly predictable, increasing process such that the Skorokhod conditions  $(\I_{\{\ygg_-\neq \Xgg_-\}}\bigcdot \lgg^r)_{\xtheta \wedge \tau}=(\I_{\{\ygg \neq \Xgg \}}\star  \lgg^g_+)_{\xtheta \wedge \tau}=0$ hold and $\lgg_0 = 0$.
\ep

\begin{proof}
The assertion of the proposition follows from Lemma \ref{lem7.1} and similar arguments as used in Section \ref{sect7} (see, in particular, the
proof of Lemma \ref{lem5.1}).
\end{proof}

We now focus on the existence of a solution to the coupled equations \eqref{eq53.0}--\eqref{eq53.1} from Proposition \ref{pro7.1}.
As in Section \ref{sect7}, we define the linear transformation
\[ \bar \xy := G\xy,\ \bar \xz:=G_{-}\xz+G^{-1}\bar \xy\nu,\
\bar \xu:=\xu,\ \bar L^r := \wt G \bigcdot L^r , \ \bar L^g := G \star L^g
\]
and the transformed generators $\bar F^{r}$ and $\bar F^{g}$. It is easy to check that if a solution $(\xy, \xz,\xu, \xl) \in \cOff \times \cPffd \times \cOff\times \overline{\mathcal{P}}(\FF)$ to the coupled equations \eqref{eq53.0}--\eqref{eq53.1} exists, then $(\bar \xy, \bar \xz, \bar \xu , \bar \xL) \in \cOff \times \cPffd \times \cOff\times \overline{\mathcal{P}}(\FF)$ satisfies the following coupled equations \eqref{eq55.0}--\eqref{eq55.1}
\begin{align}  \label{eq55.0}
\bar\xy_t&=G_\tau\xX_{\tau}-\medint\int_{\,\rrb t,\tau\rrb}\bar F_s^r(\bar\xy_s,\bar\xz_s,\bar\xu_s)\,dD^r_s-\medint\int_{\,\llb t,\tau\llb}\bar F_s^g(\bar\xy_s,\bar \xz_s,\bar \xu_s)\,dD^g_{s+} \nonumber  \\
& \quad-\medint\int_{\,\rrb t,\tau\rrb} \bar \xz_s\,d M_s+\medint\int_{\,\rrb t,\tau\rrb}\big[\wt G_s\xR_s-\Delta F^r_s\Delta D_s^r \big]d\Gamma_s-(\bar\xl_\tau-\bar\xl_t)
\end{align}
and
\begin{equation} \label{eq55.1}
\medint\int_{\,\rrb 0,t\rrb}\bar \xu_s\,d\Nogg_s=\medint\int_{\,\rrb 0,t\rrb}\big[(\xR_s-\bar \xy_sG^{-1}_s)-\wt G^{-1}_s\Delta F^r_s \Delta D_s^r\big]\,d\Nogg_s
\end{equation}
where
\[
\Delta F^r_s:= \bar F^r_s(G_s\xR_s,\bar \xz_s,\bar \xu_s)-\bar F^r_s(\bar \xy_s,\bar \xz_s,\bar \xu_s)
\]
and $\bar\xl = \bar\xl^r + \bar \xl^g$ satisfies
\[
(\I_{\{\bar\xy_- \neq G_-\xX_-\}}\bigcdot \bar \xl^r)_\tau = (\I_{\{\bar \xy \neq G\xX\}}\star \bar \xl^g)_\tau = 0.
\]
In the reverse, a solution $(\xy, \xz, \xu, \xl)$ to equations \eqref{eq53.0}--\eqref{eq53.1} can be obtained from a solution $(\bar \xy, \bar \xz ,\bar \xu)$ to equations \eqref{eq55.0}--\eqref{eq55.1} by setting $\xy := G^{-1}\bar \xy, Z := G^{-1}_-(\bar \xz - G^{-1} \bar \xy \nu)$, $\xu := \bar \xu$, $\xl^r := \wt G^{-1}\bigcdot \bar \xl^r$ and $\xl^g := G^{-1}\star \bar \xl^g$.

\section{Solution to the $\GG$ RBSDE in the Brownian case} \label{sect11}

We proceed here similarly to Section \ref{sect8}. Let $M=W$ be a $d$-dimensional Wiener process and $\FF$ be its natural filtration. We consider the coupled equations \eqref{eq53.0}--\eqref{eq53.1} with $F^g= 0$ and a c\`adl\`ag lower barrier $\xX$. The goal is again to demonstrate that, in some specific settings, the coupled equations \eqref{eq55.0}--\eqref{eq55.1} possess a unique solution $(\xy,\xz,\xu,\xl)$.

As before, we write $D$ and $F$ instead of $D^r$ and $F^r$, respectively, and we consider the case where $D=(D^1,D^2)=(\langle W \rangle,\Gamma)$
and, for all $(y,z,u) \in \RR\times \RR^d\times \RR$,
\[
F=(F^1,F^2)=(F^1({y},{z},{u}),F^2({y})).
\]
 Then the BSDE \eqref{eq53.0} is reduced to
\begin{align} \label{eq57}
\bar \xy_t& =G_\tau\xX_{\tau}-\medint\int_{\,\rrb t,\tau\rrb} \bar F_s^1\,d\langle W\rangle_s -\medint\int_{\,\rrb t,\tau\rrb} \bar F_s^2(\bar \xy_s)\,d\Gamma_s -\medint\int_{\,\rrb t,\tau\rrb} \bar \xz_s\,dW_s  \nonumber \\	
& \quad +\medint\int_{\,\rrb t,\tau\rrb}\big[\wt G_{s}\xR_s-\big(\bar F^2_s(G_s\xR_s)-\bar F^2_s(\bar \xy_s)\big)\Delta \Gamma_s \big] d\Gamma_s
-(\bar\xl_\tau-\bar\xl_t )
\end{align}
where $\bar F_s^1:=\bar F_s^1(\bar\xy_s,\bar\xz_s,\bar \xu_s)$ and equation \eqref{eq53.1} has a solution $\bar \xu$ given by equality \eqref{eq43.0}.

The following result gives sufficient conditions for existence of a solution to $\FF$ RBSDEs \eqref{eq53.0} and \eqref{eq57} in the case of the Brownian filtration $\FF$.

\bp \label{pro10.1}
Assume that:\\
(i) for every $t\geq 0$, the map $F^1_t:\RR\times \RR^{d}\times\RR\to\RR$ is bounded and Lipschitz continuous;\\
(ii) for every $t\geq 0$, the map $F^2_t:\RR\to\RR$ is bounded, Lipschitz continuous and decreasing;\\
(iii) the dual $\FF$-optional projection $A^o$ has a finite number of discontinuities;\\
(iv) the process $\xX$ is c\`adl\`ag.\\
Then the RBSDE \eqref{eq57} has a solution $(\bar\xy,\bar\xz,\bar \xl)$ and a solution $(\xy,\xz,\xl)$ to \eqref{eq53.0} can be obtained by setting $\xy:=G^{-1}\bar Y$, $\xz:= G^{-1}(\bar Z - G^{-1}\bar\xy \nu )$ and $\xl := \wt G^{-1}\bigcdot \bar\xl$.
\ep

\begin{proof}
Similarly to the proof of Proposition \ref{pro9.1},  in order to obtain a solution $(\bar\xy,\bar \xz, \bar \xl)$ to the BSDE \eqref{eq43},
we apply Theorem 2.1 in \cite{EHO2015} to the data $(\xX, R,\xtheta,F)$. We note that \eqref{eq43} is a special case of equation (2.1) in \cite{EHO2015} of the form
\begin{align*}
\bar Y_t&=G_\tau \xX_\tau+\medint\int_{\,\rrb t,\tau\rrb}f(s,\bar Y_s,\bar Z_s)\,ds+\medint\int_{\,\rrb t,\tau\rrb}g(s,\bar Y_s)\,d\bar{A}_s-\medint\int_{\,\rrb t,\tau\rrb}\bar Z_s\,dW_s
\\ & \quad +\sum_{t<s\leq\tau}h(s,\bar Y_{s-},\bar Y_s) -(\bar\xL_\tau - \bar\xl_t)
\end{align*}
where \eqref{eq57} can be recovered if we set $\bar{A}:=\Gamma^c$ (the continuous part $\Gamma$) and
\[
f(s,\bar \xy_s,\bar \xz_s):=-\bar F^1_s(\bar Y_s,\bar Z_s,\bar \xu_s), \quad g(s,\bar \xy_s):=\wt G_s\xR_s-\bar F^2_s(\bar \xy_s),
\]
and
\[
h(s,\bar \xy_{s-},\bar \xy_s):=\big(\wt G_s\xR_s-F^2_s(\bar \xy_s)\big)\Delta \Gamma_s -\big(\bar F^2_s(G_s \xR_s)-\bar F^2_s(\bar \xy_s)\big)(\Delta \Gamma_s)^2
\]
where $\bar \xu $ is given by \eqref{eq43.0}. Notice that once again $h$ does not depend on  $\bar Y_{s-}$.
As was explained in the proof of Proposition \ref{pro9.1}, the assumptions in Theorem 2.1 of \cite{EHO2015} are satisfied and thus a solution $(\bar \xy, \bar \xz, \bar \xl)$ exists. Finally, we observe that the $\FF$-predictable, increasing process $\xl := \wt G^{-1} \bigcdot \bar \xl$ clearly obeys the Skorokhod conditions since
\[
(\I_{\{\xy_- \neq \xX_-\}} \bigcdot \xl)_\tau = (\I_{\{\bar \xy_- \neq G_-\xX_-\}} \wt G^{-1} \bigcdot \bar L)_\tau = 0
\]
and thus the proof is completed.
\end{proof}

\section{BSDEs with a l\`agl\`ad driver and common jumps}   \label{sect12}

In the last part, we shall work in a general setting and study solutions to BSDEs where the driver is l\`agl\`ad and may share common jumps with the driving martingale. To the best of our knowledge, such results are not yet available in the literature.  More specifically, given a filtration $\FF$, we propose a method of solving the $\FF$ BSDE
\begin{equation} \label{eq9.0}
v_t =\xi_\tau-\medint\int_{\,\rrb t,\tau\rrb}f_s^r\,d\DB^r_s-\medint\int_{\,\llb t,\tau\llb}f_s^g\,d\DB^g_{s+}-\medint\int_{\,\rrb t,\tau\rrb}z_s\,dM_s
\end{equation}
where $f_s^r:=f_s^r(v_{s-},v_s,z_s)$ and $f_s^g:=f_s^g(v_s,v_{s+})$. In particular, we observe that in the case where either $F^g$ in \eqref{eq39.0} does not depend on $\xu$ and $\xz$ or $U$ in \eqref{eq39.1} can be solved and does not depend on $Z$ (see, for example, Section \ref{sect8}), then the BSDEs \eqref{eq39.0} and \eqref{eq40.0} can be obtained as a special case of the above BSDE \eqref{eq9.0}.

We present below a jump-adapted method of transforming the l\`agl\`ad BSDE given by \eqref{eq9.0} to a system of more tractable c\`adl\`ag BSDEs, which in turn can be further converted into a system of c\`adl\`ag BSDEs with a continuous driver. In some special cases, a solution to the latter BSDE can be obtained using results from the existing literature.
\vskip 5 pt
\noindent {\it Step 1. From a l\`agl\`ad to c\`adl\`ag driver.} For simplicity, in the following we denote $\DB:= \DB^r+ \DB^g$ so that $\DB$ is a l\`agl\`ad process of finite variation. We suppose that the times of right-hand jumps of $\DB$ (that is, the moments when $\Delta^+ \DB >0$) are given by the family $(T_i)_{i=1,2,\dots, p}$ of $\FF$-stopping times and we denote $S_0:=0,\,S_i:=T_i\wedge \tau$ and $S_{p+1}:=\tau$. We note that at each $S_i$ we have that $v_+ - v = f^g(v,v_+)\Delta \DB^g_+$, which shows that if the value of $v_+$ is already known, then the value of $v$ can be obtained as a solution to that equation. This leads to the observation that a solution $(v,z)$ to \eqref{eq9.0} can be created by first solving iteratively, for every $i = 0,1,\dots, p$, the following c\`adl\`ag BSDE on each stochastic interval $\llb S_i,S_{i+1}\rrb$
\begin{equation} \label{eq9.01}
v^i_t =\xi^i -\medint\int_{\,\rrb t,S_{i+1}\rrb}f_s^r(v^i_{s-}, v^i_s,z^i_s)\,d\DB^r_s - \medint\int_{\,\rrb t,S_{i+1} \rrb}z^i_s\,d M_s
\end{equation}
where an $\cF_{S_{i+1}}$-measurable random variable $\xi^i$ is given by the system of equations, for every $i= 0, 1,\dots, p-1$,
\[
v^{i+1}_{S_{i+1}}-\xi^i=f^g_{S_{i+1}}\big(\xi^i ,v^{i+1}_{S_{i+1}}\big)\Delta^+ \DB_{S_{i+1}}
\]
with $\xi^p=\xi_\tau$. Then a solution $(v,z)$ to the l\`agl\`ad BSDE \eqref{eq9.0} is constructed by setting
\[
v := \sum_{i=0}^p v^i\I_{\rrb S_i, S_{i+1} \rrb },\quad z := \sum_{i=0}^p z^i \I_{\rrb S_i, S_{i+1} \rrb }.
\]

\vskip 5 pt
\noindent {\it Step 2. From a c\`adl\`ag to continuous driver.}  In view of \eqref{eq9.01}, in the following we focus on showing that the c\`adl\`ag BSDE can be solved, under certain assumptions about the driver and filtration. We now consider the situation where the filtration $\FF$ can support discontinuous martingales (e.g., the Brownian-Poisson filtration) and the driver $\DB$ is possibly discontinuous. More specifically, we study the c\`adl\`ag BSDE of the form
\begin{equation} \label{eq9.3}
y_t=\xi_\tau-\medint\int_{\,\rrb t,\tau\rrb} f^r_s(y_{s-}, y_s , z)\,d\DB^r_s -\medint\int_{\,\rrb t, \tau \rrb}z_s\,dM_s
\end{equation}
where a solution $(y,z)\in \cOff\times \cPffd$ is such that $y$ is a c\`adl\`ag process.

\brem  \label{rem9.1}
Our interest in the BSDE \eqref{eq9.3} is motivated by the need to understand the well-posedness of the pre-default BSDE, which is obtained in a nonlinear reduced-form model without postulating that either condition ${\bf (C)}$ or ${\bf (A)}$ holds. The discontinuity in $\DB^r$ stems from the discontinuity of the hazard process $\Gamma$ and, in some financial applications, the introduction of the nonlinearity can be interpreted as a way to introduce ambiguity in the recovery and the default intensity (see, for instance, Fadina and Schmidt \cite{FS2019}).
\erem

In the following, we suppose that $\xi_\tau$ is bounded and $\cF_\tau$-measurable and we consider a more general BSDE
\begin{equation} \label{eq9.4}
y_t=\xi_\tau-\medint\int_{\,\rrb t,\tau\rrb}f^r_s\,d\DB^c_s -\sum_{t<s\leq \tau}h(s,y_{s-},y_s)-\medint\int_{\,\rrb t,\tau\rrb}z_s\,dM_s
\end{equation}
where $f^r_s:= f^r_s(y_{s-}, y_s,z_s)$ and $\DB^c$ is the continuous part of the process $\DB^r$.

To recover the BSDE \eqref{eq9.3} from \eqref{eq9.4}, it suffices to set $h(s,y_{s-},y_s):=f^r_s(y_{s-},y_s)\Delta \DB^r_s$. Consequently, we henceforth suppose that $h=0$ outside the graph of a finite set of $\FF$-predictable stopping times $(T_i)_{i=1,2,\dots,p}$ and we denote $S_0:=0,\,S_i:=T_i\wedge \tau$ and $S_{p+1}:=\tau$.

\brem
Note that a sufficient assumption for the jumps of $\DB^r$ to be $\FF$-predictable stopping times is to postulate that $\DB^r$ is an $\FF$-predictable,
increasing process. Furthermore, observe that the condition that $(T_i)_{i=1,2,\dots p}$ are $\FF$-predictable stopping times can be relaxed if
the mapping $h$ does not depend on $y_-$.
\erem

\brem  \label{rem9.2}
Suppose that $p=1$ and denote $S=S_1$. Let us assume that a solution to \eqref{eq9.4} on $\rrb S ,\tau \rrb$ has
already been found and our goal is to construct its extension to the interval $\llb 0,\tau \rrb$. We observe that if $(y,z)$
is a solution to \eqref{eq9.4}, then
\[
y_t=y_0+ \medint\int_{\,\rrb 0,t\rrb}f_s^r(y_{s-}, y_s,z_s)\,d\DB^c_s +\sum_{0<s\leq t}h(s,y_{s-},y_s)+\medint\int_{\,\rrb 0,t\rrb}z_s\,dM_s
\]
and hence the jump of the c\`adl\`ag process $y$ at time $S$ satisfies
\begin{equation}  \label{eq9.5}
\Delta y_S:=y_{S}-y_{S-}=h(S,y_{S-},y_S)+z_S\,\Delta M_S . 
\end{equation}

By taking the conditional expectation of both sides of \eqref{eq9.5} with respect to $\cF_{S-}$, we obtain the following equation
\[
y_{S-}=\mathbb{E}\big(y_S-h(S,y_{S-},y_{S})\,|\,\cF_{S-}\big),
\]
which, at least in principle, can be solved for $y_{S-}$ under suitable additional assumptions. Subsequently, one could compute
$z_S$ from equality \eqref{eq9.5}. However, if one decides to proceed in that way, then to solve the BSDE \eqref{eq9.4} on $\llb 0,S\rrb$
one would need to solve \eqref{eq9.4} on $\llb 0,S\llb$ and thus to study the BSDE driven by the martingale $M$ stopped at $S-$.
Since this would be quite cumbersome, we propose in Proposition \ref{pro9.2} an alternative method where this difficulty is circumvented.
\erem

To show the existence of a solution to the BSDE \eqref{eq9.4}, we introduce an auxiliary l\`agl\`ad BSDE
\begin{align} \label{eq9.6}
v_t&=\xi_\tau-h(\tau,v_{\tau-},\xi_\tau)-\medint\int_{\,\rrb t,\tau\rrb}f_s^r(v_{s-}, v_s,z_s)\, d\DB^c_s-\medint\int_{\,\rrb t,\tau\rrb}z_s\,dM_s \nonumber \\
&\quad -\sum_{t\leq s<\tau} h(s,v_{s-},v_{s+})
\end{align}
where a solution $(v,z) \in \cOff\times \cPffd$ is such that $v$ is a l\`agl\`ad process.

\bp \label{pro9.2}
Let $v$ be a l\`agl\`ad process such that $(v,z)$ is a solution to the BSDE \eqref{eq9.6} on $\llb 0,\tau\rrb$.
Then $(y,z)$ where $y:=v_+\I_{\llb 0,\tau\rrb} + h(\tau,v_{\tau-}, \xi_\tau)\I_{\llb \tau \rrb}$ is a solution to the c\`adl\`ag BSDE \eqref{eq9.4} on $\llb 0,\tau\rrb$.
\ep
	
\proof
Suppose that $(v,z)$ is a solution to \eqref{eq9.6}. It is clear from \eqref{eq9.6} that the left-hand and right-hand jumps of $v$ are given by $\Delta v=z\Delta M$ and $\Delta^+ v=h(\cdot,v_{-},v_{+})$, respectively. By the optional sampling theorem, we have that $\mathbb{E}(v_S\,|\,\cF_{S-})=v_{S-}$ for any $\FF$-predictable stopping time $S$.  Therefore, if the random variable $v_{S+}$ is known,
then the $\cF_S$-measurable random variable $v_S$ is a solution to the equation
\begin{equation}  \label{eq9.7}
\Delta^+ v_S:=v_{S+}-v_S=h\big(S,\mathbb{E}[v_S\,|\,\cF_{S-}],v_{S+}\big).
\end{equation}

If we set $y:=v_+$ on $\llb 0,\tau\llb$, then $y_{-}=v_{-}$ and thus
\begin{align*}
\Delta y_S&=y_{S}-y_{S-}=y_{S+}-y_{S-}=v_{S+}-v_{S-}=\Delta^+ v_S+\Delta v_S\\
&=h(\cdot,v_{S-},v_{S+})+z_S\,\Delta M_S=h(\cdot,y_{S-},y_{S})+z_S\,\Delta M_S,
\end{align*}
which coincides with \eqref{eq9.5}. In the next step, we take inspiration from the proof of Theorem 3.1 in Essaky et al. \cite{EHO2015} and rewrite \eqref{eq9.4} into
\begin{align*}
v_t =\xi_\tau-h(\tau,v_\tau,\xi_\tau)\Delta\DB^r_\tau-\medint\int_{\,\rrb t,\tau\rrb}f_s^r(v_{s-}, v_s,z_s)\, d\DB^c_s-\medint\int_{\,\rrb t,\tau\rrb}z_s\,dM_s
 -\sum_{t\leq s<\tau}\Delta^+ v_s.
\end{align*}
Recall that, by assumption about $h$, the right-hand jump times of $v$ are given by the family $(T_i)_{i=1,2,\dots, p}$ of $\FF$-stopping times and we denote $S_0:=0,\,S_i:=T_i\wedge \tau$ and $S_{p+1}:=\tau$. We observe that a solution $(v,z)$ can be obtained by first solving iteratively the following c\`adl\`ag BSDE on the stochastic interval $\llb S_i,S_{i+1}\rrb$, for every $i = 0,1,\dots, p$,
\begin{equation} \label{eq9.8}
v^i_t=\xi^i-\medint\int_{\,\rrb t,S_{i+1}\rrb}f^r_s(v^i_{s-}, v^i_s,z^i_s)\, d\DB^c_s - \medint\int_{\,\rrb t,S_{i+1}\rrb}z^i_s\,dM_s
\end{equation}
where $\xi^i$ is an $\cF_{S_{i+1}}$-measurable random variable determined by the recursive
system of equations, for every $i=0,1,\dots, p$ (see \eqref{eq9.7})
\begin{equation} \label{eq9.11}
v^{i+1}_{S_{i+1}}-\xi^i=h\big(S_{i+1},\mathbb{E}[\xi^i\,|\,\cF_{S_{i+1}-}],v^{i+1}_{S_{i+1}}\big)
\end{equation}
with the terminal condition $v^{p+1}_{S_{p+1}}=\xi$. In the last step, we aggregate the family of solutions $(v^i,z^i)$ for $i = 0,1,\dots, p$ by setting
\[
v:=v^0_0+\sum_{i=0}^p v^{i}\I_{\rrb S_i,S_{i+1}\rrb}, \quad z:=z^0_0+\sum_{i=0}^p z^{i}\I_{\rrb S_i,S_{i+1}\rrb}.
\]
Then, by an application of the It\^{o} formula, one can check that $(v,z)$ is a solution to the l\`agl\`ad BSDE \eqref{eq9.6}.
Furthermore, since $\Delta y_S=v_{S+}-v_{S-}$ and the dynamics of $y$ and $v$, which are given by \eqref{eq9.4} and \eqref{eq9.6},
respectively, are easily seen to coincide on each stochastic interval $\rrb S_i,S_{i+1}\llb$,
we conclude that $(y,z):=(v_+,z)$ is a solution to the c\`adl\`ag BSDE \eqref{eq9.4} once we made the appropriate adjustment to the last jump of size $h$ at the terminal time $\tau$.
\endproof

\brem \label{rem9.5}
Note that if the recovery process $R$ is $\FF$-predictable (so that one can use $A^p$ instead of $A^o$) and $D^r$ is chosen to be have $\FF$-predictable jumps (for instance, if $D^r = (\langle M \rangle, \wt G^{-1}\bigcdot A^{p}$)), then the transformed BSDE \eqref{eq40.0} has the form \eqref{eq9.4} and $h$ vanishes outside the graph of a family of $\FF$-predictable stopping times. In that case, assuming that $\xi^i$ can be solved in \eqref{eq9.11}, we would be able to consider jumps of a size $h$ depending on $v_-$.
\erem

\bex \label{ex9.2}
Let us show that if appropriate conditions are imposed on the inputs $(f^r, \DB^r, M)$, then a unique solution $(v^i,z^i)$ to \eqref{eq9.8} can be obtained on each interval $\llb S_i, S_{i+1}\rrb $ for $i =0,1,\dots, p$ and hence a solution $(y,z)$ to \eqref{eq9.4} can be constructed as well. In the following, we assume that the process $\langle M \rangle$ is continuous, the function $h$ does not depend on $v_-$ and
\[
f^r(v_-, v, z) \bigcdot \DB^c = f(v_-, v , z) \bigcdot \langle M \rangle + g(v)\bigcdot B
\]
where $B$ is an $\FF$-adapted, bounded, continuous, increasing process and $f$ and $g$ are some real-valued mappings satisfying appropriate
measurability conditions.  We note that, as $h$ does not depend on $v_-$, the assumption that the jump times of $\DB^r$ (and hence also $(S_i)_{i=1,\dots, p}$)
are $\FF$-predictable stopping times can be relaxed. Furthermore, the right-hand jumps of the process $v$ are given by $\Delta^+ v_t=h(t,v_{t+})$.

We thus need to analyze the following c\`adl\`ag BSDE with a continuous driver, on each stochastic interval $\llb S_i,S_{i+1}\rrb$ for every $i =0, 1,\dots, p$,
\begin{align*}
dv_t^{i}&=-f_t(v^i_{t-},v^i_t,z^i_t)\,d\langle M \rangle_t -g(t,v^i_t)\,dB_t-z^i_t\,dM_t,\\
v^i_{S_{i+1}}&=v^{i+1}_{S_{i+1}}-h\big(S_{i+1},v^{i+1}_{S_{i+1}}\big),
\end{align*}
with the terminal condition $v^{p+1}_{S_{p+1}}= \xi_\tau $.

\newpage

Observe that in the case of a Brownian-Poisson filtration $\FF$, the existence and uniqueness of a family of solutions $(v^i, z^i)$ can be deduced from Theorem 53.1 in Pardoux \cite{P1997} under the postulate that $f,g$ and $h$ are bounded and Lipschitz continuous functions, the process $B$ is bounded, and $M =(W,\wt N)$ where $W$ is a Brownian motion and $\wt N$ is an independent compensated Poisson process.
\eex

\bex\label{ex13.2}
Let the filtration $\FF$ be the Brownian-Poisson filtration. We consider below an example given in Gapeev et al. \cite{GJLR} of a supermartingale $J$ valued in $(0,1]$ which is the solution to the SDE
\begin{gather*}
dJ_t = -\lambda J_t\,dt + \frac{b}{\sigma} J_t (1-J_t)\,dW_t, \quad J_0 = 1.
\end{gather*}

The process $J$ takes a multiplicative form $J_t = Q_te^{-\lambda t }$ where $Q$ satisfy
\begin{gather*}
Q_t  =  1 + \medint\int^t_0 \frac{b}{\sigma}(1-J_u) Q_u\, dW_u.
\end{gather*}
For a fixed $p \in (0,1)$, we consider the supermartingales
\begin{align*}
\wt G_t& = J_t\I_{\{t \leq  T_1\}} + p J_t\I_{\{T_1 <  t\}} = J_t - (1- p)J_t\I_{\{T_1 < t\}},\\
G_t & = J_t - (1- p)J_t\I_{\{T_1 \leq t\}},
\end{align*}
and we observe that, by an application of the It\^o formula, we have
\begin{align*}
\wt G_t & = 1 + \medint\int^t_0 \frac{b}{\sigma} \wt G_{u} (1-J_u)\,dW_u  - \medint\int^t_0 \lambda \wt G_{u} \, du - (1- p)J_{T_1}\I_{\{T_1 < t\}}.
\end{align*}
We know from Jeanblanc and Li \cite{JL2020} that it is possible to construct a random time $\tau$ such that the Az\'ema optional supermartingale and the Az\'ema supermartingale associated with $\tau$ are given by $\wt G$ and $G$, respectively. In the present example, the equality $\widetilde G = G_-$ holds and the martingale $m$, the dual $\FF$-optional projection $A^o$ and the hazard process $\Gamma = (\widetilde G^{-1} \bigcdot A^o)$ associated with $\tau$ are given by the following expressions
\begin{align*}
m_t&=1 + \medint\int^t_0 \frac{b}{\sigma} \wt G_{u} (1-J_u)\, dW_u,\\
A^o_t&=\medint\int^t_0\lambda G_{u-}\, du  + (1- p)J_{T_1}\I_{\{T_1 \leq t\}},\\
\Gamma_t& =\lambda t + (1- p)\I_{\{T_1 \leq t\}},
\end{align*}
so that $\Gamma^c_t = \lambda t$ and $\Gamma^d_t = (1-p)\I_{\{T_1 \leq t\}}$.

The stopping time $T_1$ can be viewed as a shock to the underlying financial asset and $\tau$ is the timing of a default event.
The parameter $p \in (0,1)$ can be regarded as the conditional probability that the default event occurs at $T_1$ given that the default event has not occurred before $T_1$. In the following, we denote the compensated Poisson process by $\wt N$ and for the ease of presentation we set $D^r = \Gamma/(1-p)$ and  $D^g = \Gamma^d_-/(1-p) = \I_{\rrb T_1,\infty \llb}$. Furthermore, we suppose the generators $F^r$ and $F^d$ does not depend on $Y_-$ and the coupled equation \eqref{eq5.6}-\eqref{eq5.7} reduces to a single BSDE given by
\begin{align*}
\xY_t&=\xX_{T}-\medint\int_{\,\rrb t,T\rrb} \frac{F^r_s(\xY_s)}{1-p}\,d\Gamma_s-\medint\int_{\,\llb t,T\llb} \frac{F_s^g (\xY_s)}{1-p}\,d\Gamma^d_{s} + \medint\int_{\,\rrb t,T\rrb} \frac{b}{\sigma}(1-J_s) Z^1_s\, ds -\medint\int_{\,\rrb t,T\rrb}Z^1_s\,dW_s \\
& \quad -\medint\int_{\,\rrb t, T\rrb}Z^2_s\,d\wt N_s  +\medint\int_{\,\rrb t,T\rrb}\big[\xR_s-\xY_s-(F^r_s(\xR_s)-F^r_s(\xY_s))\I_{\llb T_1\rrb}(s)\big] d\Gamma_s.
\end{align*}

On the set $\{T_1 \leq T\}$, we observe that the driver of the above BSDE has only one jump at time $T_1$ and thus on the stochastic interval $\rrb T_1 , T\rrb$ we need only to find the solution $(y,u)$ where $u = (u^1,u^2)$ to the BSDE,
\begin{align} \label{yT1}
y_t    & = X_T - \medint\int^T_t \Big[\frac{\lambda F^r_s(y_s)}{1-p} - \frac{b(1-J_s)u^1_s }{\sigma}- \frac{\lambda(R_s - y_s)}{1-p}\Big]\,ds\nonumber \\
	   & \quad +  \medint\int^T_t u^1_s\, dW_s + \medint\int^T_t u^2_s\, d\wt N_s.
\end{align}
At the jump time $T_1$, the right jump of $Y$ is given by $\Delta^+Y_{T_1} = F_{T_1}^g(Y_{T_1})$ and the quantity $Y_{T_1}$ is obtained by solving the equation $y_{T_1} - F_{T_1}^g (Y_{T_1}) = Y_{T_1}$.

Given that $Y_{T_1}$ can be obtained, we see that one is required to solve the c\`adl\`ag BSDE, on the stochastic  interval $\llb 0, T_1\rrb$,
\begin{align*}
\xY_t & =Y_{T_1}-\medint\int_{\,\rrb t,T_1\rrb} \frac{F^r_s(\xY_s)}{1-p}\,d\Gamma_s + \medint\int_{\,\rrb t,T_1\rrb}\frac{b}{\sigma}(1-G_s)\,Z^1_s ds-\medint\int_{\,\rrb t,T_1\rrb} Z^1_s\,dW_s -\medint\int_{\,\rrb t, T_1\rrb} Z^2_s\,d\wt N_s  \\
&\quad +\medint\int_{\,\rrb t,T_1\rrb}\big[\xR_s-\xY_s-(F^r_s(\xR_s)-F^r_s(\xY_s))\I_{\llb T_1 \rrb}(s)\big]\, d\Gamma_s,
\end{align*}
which is a GBSDE where the martingale term and the driver can share a common jump at $T_1$. Again, we observe that the driver jumps only at $T_1$ with the jump size given by
\begin{gather*}
h(T_1,Y_{T_1}): = F^r_{T_1}(Y_{T_1}) - \big[\xR_{T_1}-Y_{T_1}-(F^r_{T_1}(\xR_{T_1})-F^r_{T_1}(Y_{T_1}))\big](1-p).
\end{gather*}
Therefore, the adjusted terminal condition $v_{T_1}$ at $T_1$ equals
\begin{align*}
v_{T_1} &:= Y_{T_1} - h(T_1, Y_{T_1}) = \xR_{T_1}-F^r_{T_1}(\xR_{T_1}) - p\big[\xR_{T_1}-Y_{T_1}-(F^r_{T_1}(\xR_{T_1})-F^r_{T_1}(Y_{T_1})\big]
\end{align*}
and we see that we need to solve the following BSDE with a continuous driver, on the stochastic interval $\llb 0,  T_1 \rrb$,
\begin{align}\label{vT1}
v_t		& = v_{T_1} -\medint\int_{\,\rrb t,T_1\rrb} \Big[ \frac{\lambda  F^r_s(v_s)}{1-p}+ \frac{b(1-J_s)\,z^1_s }{\sigma}+ \frac{\lambda(\xR_s-v_s)}{1-p}\Big]\,ds \nonumber\\
		& \quad -\medint\int_{\,\rrb t,T_1\rrb} z^1_s\,dW_s -\medint\int_{\,\rrb t, T_1\rrb} z^2_s\,d\wt N_s.
\end{align}
To this end, let $Y_{T_1}$ be a solution to the equation $y_{T_1} - F_{T_1}^g (Y_{T_1}) = Y_{T_1}$.
Then a solution $(Y,Z)$ where $Z = (Z^1, Z^2)$ on the whole interval $\llb 0 ,T\rrb$ can be obtained by setting
\begin{align*}
Y&:= v\I_{\llb 0, T_1\rrb}+ h(T_1, Y_{T_1})\I_{\llb T_1 \rrb } + y\I_{\rrb T_1, T\rrb},\\
Z^i&:= z^i\I_{\llb 0, T_1\rrb} + u^i\I_{\rrb T_1 ,  T\rrb}.
\end{align*}
Let us now consider the set $\{T_1 > T\}$. Since there are no jumps before $T$, it suffices to find $(v,z)$ in \eqref{vT1} on the whole interval $\llb 0 ,T\rrb$ with the terminal condition $v_T = X_T$.

To better visualise the jump-adapted method outlined above, we give a graphical illustration
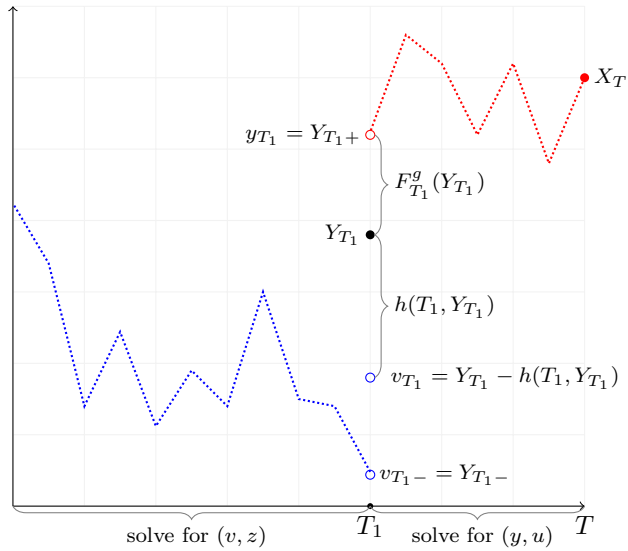
\begin{center}
\begin{tikzpicture}[scale=1.9]
\draw [step=0.5,thin,gray!10] (0,0) grid (4,3.5);

\draw [->] (0,0) -- (4,0) node [below] {$T$};
\draw [->] (0,0) -- (0,3.5) node [left] {};

\draw [densely dotted, red, thick] (4,3) -- (3.75,2.4) -- (3.5, 3.1) -- (3.25,2.6) -- (3,3.1) -- (2.75,3.3) -- (2.5,2.62);
\node (a) at (2.5,2.6) [draw,circle, red, inner sep=1.2pt]{} ;
\draw (a) node [left] {\footnotesize $y_{T_1}= Y_{T_1+}$} ;
\node (b) at (2.5,1.9) [fill, circle,inner sep=1.2pt]{} ;
\draw (b) node [left] {\footnotesize $Y_{T_1}$} ;
\node (c) at (2.5,0.9) [draw,blue,circle,inner sep=1.2pt]{} ;
\draw (c) node [right] {\footnotesize $\,\,\,v_{T_1} = Y_{T_1} - h(T_1, Y_{T_1})$} ;
\node (d) at (4,3) [fill,circle,red,inner sep=1.2pt]{} ;
\draw (d) node [right] {\footnotesize $X_{T}$} ;
\node (e) at (2.5,0) [fill,circle,inner sep=0.8pt]{} ;
\draw (e) node [below] {$T_1$} ;
\node (f) at (2.5,0.22) [draw, blue, circle,inner sep=1.2pt]{} ;
\draw (f) node [right] {\footnotesize  $v_{T_1-} = Y_{T_1-}$} ;
\draw [densely dotted, blue,thick] (2.5,0.24) -- (2.25, 0.7)  -- (2,0.75)   -- (1.75,1.5)  -- (1.5,0.7)  -- (1.25,0.95)  -- (1,0.56)  -- (0.75,1.22) -- (0.5,0.7)-- (0.25,1.7)-- (0,2.12) ;

\draw [gray,decorate,decoration={brace,amplitude=5pt}]
   (2.53,2.6) -- (2.53,1.9)
   node [black,midway,right=2pt,xshift=2pt]
   {\footnotesize $F^g_{T_1}(Y_{T_1})$};

   \draw [gray,decorate,decoration={brace,amplitude=5pt}]
   (2.53,1.9) -- (2.53,0.9)
   node [black,midway,right=2pt, xshift=2pt]
   {\footnotesize $h(T_1,Y_{T_1})$};

   \draw [gray,decorate,decoration={brace,amplitude=5pt}]
   (4,0) -- (2.5,0)
   node [black,midway,below=4pt, xshift=2pt]
   {\footnotesize solve for $(y,u)$};

      \draw [gray,decorate,decoration={brace,amplitude=5pt}]
   (2.5,0) -- (0,0)
   node [black,midway,below=4pt, xshift=2pt]
   {\footnotesize solve for $(v,z)$};
\end{tikzpicture}
\end{center}
Here we point out that since $(1-J)$ is bounded by one, when considering \eqref{yT1} and \eqref{vT1}, we do not need to study the transformed BSDE given in \eqref{eq40.0}-\eqref{eq40.1}. This is because, given appropriate assumptions on $F^r$, the linear growth conditions in $z$ can be easily verified here.

\eex
\section{RBSDEs with a l\`agl\`ad driver and common jumps} \label{sect13}

Following the structure of Section \ref{sect12}, given a filtration $\FF$, we focus on $\FF$ RBSDEs of the form
\begin{align}
v_t &=  \xi_\tau - \medint\int_{\,\rrb t,\tau\rrb} f_s^r(v_{s-}, v_s,z_s)\,d\DB^r_s - \medint\int_{\,\llb t, \tau \llb} f_s^g(v_s,v_{s+})\,d\DB^g_{s+} \nonumber
\\&\quad  -\medint\int_{\,\rrb t,\tau\rrb} z_s\,dM_s + l^r_\tau -l^r_{t}+l^g_\tau-l^g_{t} \label{eq10.0}
\end{align}
where $l^r$ and $l^g$ satisfy $(\I_{\{v_-\neq \xi_-\}}\bigcdot l^r)_\tau=(\I_{\{v \neq \xi\}}\star l^g)_\tau=0$.
We observe that in the case where $F^g$ in \eqref{eq53.0} does not depend on $\xu$ and $\xz$ or that $U$ in \eqref{eq53.1} can be solved and does not depend on $Z$ (for an example, see Section \ref{sect11}), then both BSDE \eqref{eq53.0} and \eqref{eq55.0} can be obtained as a special case of the above BSDE \eqref{eq10.0}. Similar to the non-reflected case, we present below a jump-adapted method to reduce the l\`agl\`ad RBSDE \eqref{eq10.0} to a system of c\`adl\`ag RBSDEs, which can be further reduced to a system of c\`adl\`ag RBSDEs with continuous drivers.

\vskip 5 pt
\noindent {\it Step 1. From a l\`agl\`ad to c\`adl\`ag driver.}  By examining the right-hand jumps of $v$, that is, $\Delta^+v$, and the Skorokhod condition satisfied by $l^g$, we observe that $\Delta^+ v$ and $\Delta l^g_+$ must satisfy the conditions
\[
v_+ - v  = f^g(v_+, v)\Delta \DB^g_+ + \Delta l^g_+ , \quad (v - \xi)\Delta l^g_+=0,
\]
which in turn implies that
\[
\Delta l^g_+=\left(\xi-(v_+ - f^g(v,v_+)\Delta \DB^g_+)\right)^+ , \quad v=\xi\vee \left(v_+ - f^g(v, v_+)\Delta \DB^g_+\right).
\]
We thus see that, at the jump times of $l^g_+$ and $\DB^g_+$, the quantity $v$ can be obtained by solving the second equation (of course,
assuming that a solution exists) and $\Delta l^g_+$ can be obtained by substitution. In particular, if $f^g$ does not depend on $v$, then
it is clear that we have
\[
\Delta l^g_+ =(\xi-(v_+ -f^g(v_+)\Delta \DB^g_+) )^+, \quad v=\xi\vee (v_+ -f^g(v_+)\Delta \DB^g_+).
\]
These arguments lead to the observation that a solution $(v,z,l )$ to \eqref{eq10.0} can be obtained by solving iteratively, for every $i = 0,1,\dots, p$, the following c\`adl\`ag RBSDE on $\llb S_i,S_{i+1}\rrb$
\begin{equation} \label{eq10.01}
v^i_t=\xi^i-\medint\int_{\,\rrb t,S_{i+1}\rrb}h_s^r(v^i_{s-},v^i_s,z^i_s)\,d\DB^r_s-\medint\int_{\,\rrb t,S_{i+1}\rrb}z_s^i\,dM_s+l^i_{S_{i+1}}-l^i_{t}
\end{equation}
where the c\`adl\`ag increasing process $l^i$ obeys the Skorokhod condition $(\I_{\{\xi_- =v^i_-\}}\bigcdot l^i)=0$
and  $(\xi^i, \Delta l^g_{S_{i+1}+})$ are $\cF_{S_{i+1}}$-measurable random variables such that $\Delta l^g_{S_{p+1}+} = 0$, $\xi^p =\xi_{\tau}$ and for $i = 0,1,\dots, p-1$,
\begin{align*}
\Delta l^g_{S_{i+1}+}&=\left( \xi_{S_{i+1}}-(v_{S_{i+1}}^{i+1}-h^g_{S_{i+1}}(\xi^i,v_{S_{i+1}}^{i+1})\Delta \DB^g_{S_{i+1}+})\right)^+,\\
\xi^{i}&=\xi_{S_{i+1}}\vee \left(v^{i+1}_{S_{i+1}}-h^g_{S_{i+1}}(\xi^i,v^{i+1}_{S_{i+1}})\Delta \DB^g_{S_{{i+1}}+}\right).
\end{align*}
Then a global solution $(v,z,l)$ where $l = l^r + l^g$ is obtained by setting
\begin{align*}
v & =v_0+\sum_{i=0}^p v^i\I_{\rrb S_i,S_{i+1}\rrb },\qquad z=z_0+\sum_{i=0}^p z^i\I_{\rrb S_i,S_{i+1} \rrb },\\
l^r & =\sum_{i=0}^p(l^{i-1}_{S_i}+l^i)\I_{\rrb S_i,S_{i+1} \rrb },\quad l^g=\sum_{i=1}^p\Delta l^g_{S_{i}+}\I_{\rrb S_{i},\infty \llb},
\end{align*}
where $l^{-1}_0 =0,\, v_0 = v^0_0$ and $z_0 = z^0_0$.

\vskip5pt
\noindent {\it Step 2. From a c\`adl\`ag to continuous driver.} In view of the c\`adl\`ag RBSDE \eqref{eq10.01}, we study the RBSDE of the form
\begin{equation} \label{eq10.3}
y_t=\xi_{\tau }-\medint\int_{\,\rrb t,\tau \rrb}f^r_s(y_{s-}, y_s,z_s)\,d\DB^r_s -\medint\int_{\,\rrb t, \tau \rrb}z_s\,dM_s + l_{\tau} - l_t
\end{equation}
where a solution $(y,z,l)\in \cOff\times \cPffd\times \cPff$ is such that $y$ is a c\`adl\`ag process and $l$ is a c\`adl\`ag, increasing process  such that $(\I_{\{y_-\neq \xi_-\}}\bigcdot l)_\tau=0$ and $l_0=0$. In the following, we consider a more general RBSDE of the form
\begin{align} \label{eq10.4}
y_t=\xi_\tau-\medint\int_{\,\rrb t,\tau\rrb}f_s^r(y_{s-}, y_s,z_s)\,d\DB^c_s-\sum_{t<s\leq \tau} h(s,y_{s-},y_s) -\medint\int_{\,\rrb t,\tau\rrb}z_s\,dM_s + l_\tau- l_t
\end{align}
where a solution $(y,z,l)\in \cOff\times \cPffd\times \cPff$ is such that $y$ is a c\`adl\`ag and $l$ is a c\`adl\`ag increasing process such that
$(\I_{\{y_- \neq \xi_-\}}\bigcdot l^r)_\tau=0$ and $l_0 = 0$.

To recover equation \eqref{eq10.3} from \eqref{eq10.4}, it suffices to set $h(s,y_{s-},y_s):=f^r_s(y_{s-}, y_s)\Delta\DB_s^r$. In view of this, we further suppose that $h=0$ outside the graph of a finite family of $\FF$-predictable stopping times $(T_i)_{i=1,2,\dots,p}$ and we denote $S_0=0,\,S_i=T_i\wedge \tau$ and $S_{p+1}=\tau$. To examine the existence of a solution to the RBSDE \eqref{eq10.4}, we introduce an auxiliary RBSDE
\begin{align} \label{eq10.5}
v_t&=\xi_\tau-h(\tau,v_{\tau-},\xi_\tau)-\medint\int_{\,\rrb t,\tau\rrb}f_s^r(v_{s-}, v_s,z_s)\,d\DB^c_s  \nonumber \\
&\quad -\sum_{t\leq s<\tau} h(s,v_{s-},v_{s+})-\medint\int_{\,\rrb t,\tau\rrb}z_s\,dM_s + l_\tau - l_t
\end{align}
where a solution $(v,z,l)$ is such that $v$ is a l\`agl\`ad, $\FF$-adapted process, the process $z$ is $\FF$-predictable and the process $l$ obeys the Skorokhod condition $(\I_{\{v_- \neq \xi_-\}}\bigcdot l)_\tau=0$ and $l_0 = 0$.

\bp \label{pro10.2}
Let $v$ be a l\`agl\`ad process such that $(v,z,l)$ is a solution to the RBSDE \eqref{eq10.5} on $\llb 0,\tau\rrb$.
Then $(y,z,l)$ where $y:=v_+\I_{\llb 0,\tau\rrb} + h(\tau, v_{\tau-} , \xi_{\tau})\I_{\llb \tau \rrb}$ solves the c\`adl\`ag RBSDE \eqref{eq10.4} on $\llb 0,\tau\rrb$.
\ep

\proof
Suppose that $(v,z,l)$ is a solution to \eqref{eq10.4}. It is clear from \eqref{eq10.4} that the left-hand and right-hand jumps
of $v$ are given by $\Delta v=z\Delta M + \Delta l$ and $\Delta^+ v=h(\cdot,v_{-},v_{+})$, respectively.
Note that $l$ must satisfy the reflection condition $(v_{S-} - \xi_{S-})\Delta l_S = 0$ and, by the optional sampling theorem, we have that $\mathbb{E}[v_S\,|\,\cF_{S-}]=v_{S-} + \Delta l_S$ for any $\FF$-stopping time $S$. Then, by solving these two equations, we obtain
\[
v_{S-}=\xi_{S-}\vee \mathbb{E}[v_S\,|\,\cF_{S-}] , \quad  \Delta l_S =\big(\xi_{S-}-\mathbb{E}[v_S \,|\,\cF_{S-}]\big)^+.
\]
Therefore, if the random variable $v_{S+}$ is known, then the $\cF_S$-measurable random variable $v_S$ is a solution to the equation
\[
\Delta^+ v_S:=v_{S+}-v_S=h\big(S,\xi_{S-}\vee \mathbb{E}[v_S \,|\,\cF_{S-}],v_{S+}\big).
\]
Recall that, by assumption about $h$, the right-hand jump times of $v$ are given by the family $(T_i)_{i=1,2,\dots, p}$ of $\FF$-predictable stopping times and we set $S_0=0,\,S_i=T_i\wedge \tau$ and $S_{p+1}=\tau$. We observe that a solution $(v,z,l)$ can be constructed by first solving by iteration, for every $i = 0,1,\dots, p$, the following c\`adl\`ag RBSDE on the stochastic interval $\llb S_i,S_{i+1}\rrb$
\[
v^i_t=\xi^i-\medint\int_{\,\rrb t,S_{i+1}\rrb}f_s^r(v^i_{s-}, v^i_s,z^i_s)\,d\DB^c_s -\medint\int_{\,\rrb t,S_{i+1}\rrb}z^i_s\,dM_s + l^i_{S_{i+1}} - l^i_t
\]
where, for every $i=0,1,\dots, p$, we have $\I_{\{v^i_- \neq \xi_- \}\cap \rrb S_i,S_{i+1}\rrb}\bigcdot l^i=0$ and we denote by $\xi^i$ an $\cF_{S_{i+1}}$-measurable random variable satisfying  $v^{p+1}_{S_{p+1}}=\xi_\tau$ and
\[
v^{i+1}_{S_{i+1}}-\xi^i=h\big(S_{i+1},\xi_{S_{i+1}-}\vee\mathbb{E}[\xi^i\,|\,\cF_{S_{i+1}-}],v^{i+1}_{S_{i+1}}\big).
\]

We aggregate the family of solutions $(v^i,z^i)$ for $i = 0,1,\dots, p$ by setting $l^{-1}_0 = 0$ and
\[
v=v^0_0+\sum_{i=0}^p v^{i}\I_{\rrb S_i,S_{i+1}\rrb},\quad z=z^0_0+\sum_{i=0}^p z^{i}\I_{\rrb S_i,S_{i+1}\rrb}, \quad l=\sum_{i=0}^p (l^{i-1}_{S_i} + l^{i})\I_{\rrb S_i,S_{i+1}\rrb}.
\]
Using the It\^{o} formula, one can check that $(v,z,l)$ satisfies the l\`agl\`ad RBSDE \eqref{eq10.4} and $l$ is increasing and satisfies $(\I_{\{v_- \neq \xi_-\}}\bigcdot l)_\tau=0$. Furthermore, since $\Delta y_S=v_{S+}-v_{S-}$ and the dynamics of $y$ and $v$
(see \eqref{eq10.4} and \eqref{eq10.5}, respectively) are easily seen to coincide on $\rrb S_i,S_{i+1}\llb$, we conclude that $(y,z,l):=(v_+,z,l)$ is a solution to the c\`adl\`ag RBSDE \eqref{eq10.4} once we made the appropriate adjustment to the last jump of size $h$ at the terminal time $\tau$, since $v_-=y_-$ and $(\I_{\{v_-\neq \xi_-\}}\bigcdot l)_\tau=0$.
\endproof

\bex \label{ex10.1}
Here we show that if appropriate conditions are imposed on the inputs data $(f,\DB, M)$, then a unique solution $(v^i,z^i)$ to \eqref{eq9.8} can be obtained on each interval $\llb S_i, S_{i+1}\rrb $ for $i =0,1,\dots, p$ and thus a solution $(y,z,l)$ to \eqref{eq9.4} can be constructed. In the following, we assume that the process $\langle M \rangle$ is continuous, the function $h$ does not depend on $v_-$ and
\[
f^r(v_-, v, z) \bigcdot \DB^c = f(v_-, v , z) \bigcdot \langle M \rangle + g(v)\bigcdot C
\]
where $C$ is an $\FF$-adapted, continuous, increasing process. Furthermore, $f$ and $g$ are some real-valued mappings that satisfies appropriate measurability conditions.

We note that since $h$ does not depend on $v_-$, the assumption that the jumps of $\DB$ occur at $\FF$-predictable stopping times can be relaxed
and the right-hand jumps of $v$ are given by $\Delta^+ v_t=h(t,v_{t+})$. Hence one is required to solve the following c\`adl\`ag RBSDE with continuous drivers, on each stochastic interval $\llb S_i,S_{i+1}\rrb$ for $i = 0,1,\dots, p$,
\begin{align*}
dv_t^{i}&=-f_t(v^i_t,z^i_t)\,dt-g_t(v^i_t)\,d\DB_t^c-z^i_t\,dM_t + dl^i_t,\\
v^i_{S_{i+1}}&=v^{i+1}_{S_{i+1}}-h\big(S_{i+1},v^{i+1}_{S_{i+1}}\big),
\end{align*}
where $v^{p+1}_{S_{p+1}}=\xi $ and $l^i \in \cPff$ is a c\`adl\`ag, increasing process with $l^i_0=0$ and such that the following equality holds
\[
\I_{\{v^i_-\neq \xi_-\}\cap \rrb S_i,S_{i+1}\rrb}\bigcdot l^i=0.
\]
In the case where $\FF$ is a Poisson filtration or, more generally, is generated by the {\it Teugels martingales} (see Nualart and Schoutens \cite{NS2000}
or Schoutens and Teugels \cite{ST1998}), the existence and uniqueness of a solution $(v^i, z^i,l^i)$ can be obtained by an application of Theorem 5 in Ren and El Otmani \cite{RE2010} under the postulate that $f, g$ and $h$ are bounded and Lipschitz continuous functions, the process $\DB$ is bounded and $M$ is the compensated Poisson process.
\eex

\section{Appendix}       \label{sect14}

We assume that the process $R$ is $\FF$-optional and we define the l\`agl\`ad process $Q$ by
\begin{align*}
Q :=\xN-\xR\bigcdot A^o+C=\xN-\xR\bigcdot A^o+C^r+C^g
\end{align*}
where $\xN$ is an $\FF$-local martingale and $C$ is a l\`adl\`ag process of finite variation.
If $Y$ is a l\`adl\`ag process of finite variation or, more generally, an optional semimartingale
(which, by definition, is assumed to be a l\`adl\`ag process), then $Y$ admits the decomposition $Y=Y^r+Y^g$
where $Y^g_t:=\sum_{s<t}(Y_{s+}-Y_s)$ and the c\`adl\`ag process $Y^r$ is given by $Y^r:=Y-Y^g$.

\bl \label{lem11.1}
Assume that $G>0$. Then the process $G^{-1}$ satisfies
\begin{equation} \label{eq11.1}
G^{-1}=G_0 -G_{-}^{-2}\bigcdot \wt m+G^{-1}\bigcdot \Gamma
\end{equation}
where $\Gamma := \wtG^{-1} \bigcdot A^o$ and $\wt m:=m-\wtG^{-1}\bigcdot [m,m]$. Moreover, for the c\`adl\`ag process
\[
Q^r := \xN-\xR\bigcdot A^o+C^r
\]
we have that
\begin{equation} \label{eq11.2}
[Q^r, G^{-1}]=-G^{-1}G^{-1}_{-}\big([\xN,m]-[\xN,A^o]+[C^r,G]\big)-R\Delta G^{-1}\bigcdot A^o.
\end{equation}
\el

\begin{proof}
For brevity, we write $[G]:=[G,G]$ and $[m]=[m,m]$. The It\^o formula yields
\begin{equation} \label{eq11.2a}
G^{-1}=G^{-1}_0 -G_{-}^{-2}\bigcdot G+G^{-1}G^{-2}_{-}\bigcdot [G]= G^{-1}_0 -G_{-}^{-2}\bigcdot J
\end{equation}
where $J:=G-G^{-1}\bigcdot [G]$. Since $G=m-A^o$ and thus $\Delta G=\Delta m-\Delta A^o$, we obtain
\begin{align*}
[G]& =[m]-[m,A^o]+[A^o,A^o]=[m]-\Delta m\bigcdot A^o-(\Delta m-\Delta A^o)\bigcdot A^o\\
&=[m]-\Delta m\bigcdot A^o-\Delta G\bigcdot A^o
\end{align*}
so that
\begin{align*}
J=G-G^{-1}\bigcdot [G]=m_-A^o-G^{-1} \bigcdot [m]+G^{-1}\big(\Delta m\bigcdot A^o + \Delta G\bigcdot A^o\big).
\end{align*}
Using the equalities $\wt m =m-\wtG^{-1}\bigcdot [m]$ and $\Delta m=\wtG-G_{-}$, we get
\begin{align*}
J= \wt{m}+\wtG^{-1}\bigcdot [m]-A^o-G^{-1}\bigcdot [m]- G^{-1}\big((\wtG-G_{-})\bigcdot A^o-\Delta G\bigcdot A^o\big).
\end{align*}
Since $\wtG-G=\Delta A^o$, we also have that
\begin{align*}
(\wtG^{-1}-G^{-1})\bigcdot [m]&=G^{-1}\wtG^{-1}(G-\wtG)\bigcdot [m]=-G^{-1}\wtG^{-1}\Delta A^o\bigcdot [m] \\
&=-G^{-1}\wtG^{-1}(\Delta m)^2\bigcdot A^o=-G^{-1}\wtG^{-1}(\wtG-G_{-})^2\bigcdot A^o.
\end{align*}
Consequently,
\begin{align*}
J=\wt{m}-A^o+G^{-1}\big(\wtG-G_{-}+\Delta G-\wtG^{-1}(\wtG-G_{-})^2\big)\bigcdot A^o =\wt{m}-G^{-1} G^2_{-}\bigcdot \Gamma_,
\end{align*}
which, when combined with \eqref{eq11.2a}, shows that \eqref{eq11.1} is valid. To establish \eqref{eq11.2}, we first compute
\begin{align*}
[Q^r, G^{-1}]&=-G^{-2}_{-}\bigcdot [Q^r,G]+G^{-1}G_{-}^{-2}\bigcdot [Q^r,[G]]\\
&=-G^{-2}_{-}\bigcdot [Q^r,G]+G^{-1}G_{-}^{-1}\Delta G\bigcdot [Q^r,G]=-G^{-1}G^{-1}_{-}\bigcdot [Q^r,G].
\end{align*}
Finally, using the equalities $\Delta G=-GG_{-}\Delta G^{-1}$ and $G=m-A^o$, we obtain
\begin{align*}
[Q^r, G^{-1}]&=-G^{-1}G_{-}^{-1}\bigcdot [Q^r,G]=G^{-1}G_{-}^{-1}R\bigcdot [A^o,G] -G^{-1}G_{-}^{-1}\bigcdot \big([\xN,G]+[C^r,G]\big) \\
&=-G^{-1}G_{-}^{-1}\bigcdot \big([\xN,G]+[C^r,G]\big)-R\Delta G^{-1}\bigcdot A^o\\
&=-G^{-1}G_{-}^{-1}\bigcdot ([\xN,m]-[\xN,A^o]+[C^r,G])-R\Delta G^{-1}\bigcdot A^o
\end{align*}
and thus equality \eqref{eq11.2} is proven as well.
\end{proof}

We maintain the assumption that $G>0$ (and thus also $\wtG>0$) and we consider the process
\begin{equation} \label{eq11.3}
Y:=G^{-1}Q=G^{-1}(\xN-\xR\bigcdot A^o+C).
\end{equation}
Our goal is to derive the dynamics of $Y$ in terms of $\Gamma, \wt m$ and
\begin{align*}
\wt\xN:=\xN-\wtG^{-1}\bigcdot [\xN,m].
\end{align*}
In the proof of Lemma \ref{lem11.3}, we will employ the optional integration by parts formula. Recall that,
by definition, any {\it semimartingale} is a c\`adl\`ag process but an {\it optional semimartingale} is not necessarily a c\`adl\`ag process
although, by definition, it is a l\`adl\`ag process.

Let $X=X^r+X^g$ and $Y=Y^r+Y^g$ be l\`agl\`ad optional semimartingales such as $Y$ is of finite variation.
Then the {\it optional integration by parts formula} reads (see Theorem 8.2 in Gal'\v{c}uk \cite{G1981})
\begin{equation} \label{eq11.4}
XY= X_0Y_0+ X \circ Y + Y\circ X + [X,Y]
\end{equation}
where the {\it optional stochastic integrals} are given by
\begin{align*}
&X\circ Y=X_{-}\bigcdot Y^r+X \star Y^g_{+},\\
&Y\circ X=Y_{-}\bigcdot X^r+Y \star X^g_{+},
\end{align*}
where $X^g_+$ (respectively, $Y^g_+$) is the c\`adl\`ag version of the c\`agl\`ad process $X^g$ (respectively, $Y^g$)
and the quadratic covariation $[X,Y]$ equals
\[
[X,Y]_t=\sum_{0<s\leq t}\Delta X_s\Delta Y_s+\sum_{0\leq s<t}\Delta^+ X_s\Delta^+Y_s
\]
where we denote $\Delta X_t=X_t-X_{t-}$ and $\Delta^+ X_t =X_{t+}-X_{t}$.

For the reader's convenience, we formulate a variant of the optional integration
by parts formula \eqref{eq11.4}, which holds when $X=X^r$ is a (c\`adl\`ag) semimartingale and $Y=Y^g$ is a c\`agl\`ad process of finite variation.

\bl \label{lem11.2}
Let $X=X^r$ be a semimartingale and let $Y=Y^g$ be a c\`agl\`ad process of finite variation.  Then the process
$XY$ is l\`adl\`ag and satisfies, for every $0 \leq s< t$,
\begin{equation} \label{eq11.5}
X_tY_t=X_sY_s+\medint\int_{\,\rrb s,t \rrb}Y_{u}\,dX_u+\medint\int_{\,\llb s,t\llb }X_u\,dY^g_{u+}
\end{equation}
where $Y^g_+$ is the c\`adl\`ag version of $Y^g$.
\el

We will use the shorthand notation for \eqref{eq11.5}
\begin{equation*}
XY= X_0Y_0+ Y \bigcdot X +X \star Y^g_{+}
\end{equation*}
but all equalities in the proof of Lemma \ref{lem11.3} should be understood in the sense of \eqref{eq11.5}, meaning that all integrals with
respect to a c\`adl\`ag  (respectively, c\`agl\`ad) process should be evaluated on the interval $\rrb s,t\rrb$ (respectively, on the interval $\llb s,t\llb$) for arbitrary $0 \leq s< t$.

\bl \label{lem11.3}
If the process $Y$ is given by \eqref{eq11.3} where $C$ is a l\`adl\`ag process of finite variation with the decomposition $C=C^r+C^g$, then
\begin{align} \label{eq11.6}
Y_t&=Y_0-\medint\int_{\,\rrb 0,t\rrb}(\xR_s-Y_s)\,d\Gamma_s-\medint\int_{\,\rrb 0,t\rrb}Y_{s-}G_{s-}^{-1}\,d\wt m_s+\medint\int_{\,\rrb 0,t\rrb } G^{-1}_{s-}\,d\wt\xN_s\\
&+ \medint\int_{\,\rrb 0,t\rrb }\wtG^{-1}_s\,dC^r_s+\medint\int_{\,\llb 0,t\llb }G^{-1}_s\,dC^g_{s+}. \nonumber
\end{align}
\el

\begin{proof}
We note that $Q$ satisfies
\begin{align*}
Q=\xN-\xR\bigcdot A^o+C=\xN-\xR\bigcdot A^o+C^r+C^g=Q^r+C^g
\end{align*}
where
\begin{align*}
Q^r=\xN-\xR \bigcdot A^o+C^r =\wt \xN +\wtG^{-1}\bigcdot [\xN,m]-\xR\bigcdot A^o+C^r .
\end{align*}
The integration by parts formulas applied to $Y=G^{-1}Q=G^{-1}Q^r+G^{-1}C^g$ gives
\begin{align} \label{eq11.7}
Y&=G^{-1}Q^r+G^{-1}C^g \nonumber \\
&= Y_0 + Q_{-}\bigcdot G^{-1}+G^{-1}_{-}\bigcdot Q^r+[Q^r,G^{-1}]+G^{-1}\star C^g_{+}
\end{align}
since $G^{-1}$ and $Q^r$ are (c\`adl\`ag) semimartingales and thus the It\^o integration by parts formula
applied to the product $G^{-1}Q^r$ yields
\[
G^{-1}Q^r=Q^r_{-}\bigcdot G^{-1}+G^{-1}_{-}\bigcdot Q^r+[Q^r,G^{-1}]
\]
whereas the optional integration by parts formula (see Lemma \ref{lem11.2}) gives
\[
G^{-1}C^g=C^{g}\bigcdot G^{-1}+G^{-1} \star C^g_{+}.
\]
From \eqref{eq11.2} and \eqref{eq11.7}, we obtain
\begin{align} \label{eq11.8}
Y&=Y_0-\,Q_{-}\big(G_{-}^{-2}\bigcdot \wt m-G^{-1}\bigcdot \Gamma \big) +G^{-1}_{-}\bigcdot \big(\wt \xN +\wtG^{-1}\bigcdot [\xN,m]-\xR\bigcdot A^o+C^r\big) \nonumber \\ &\quad \, -G^{-1}_{-}G^{-1}\big( [\xN,m]- [\xN,A^o]+ [C^r,G]\big) -R\Delta G^{-1}\bigcdot A^o +G^{-1}\star C^g_{+} \\
&=Y_0-Y_{-}G_{-}^{-1}\bigcdot \wt m+G^{-1}_{-}\bigcdot \wt\xN+ K+ H \nonumber
\end{align}
where
\begin{align*}
K &:=G^{-1}_{-}\bigcdot C^r + G^{-1} \star C^g_{+} - G^{-1}G^{-1}_{-}\bigcdot [C^r,G] \\
& = G^{-1}_{-}\bigcdot C^r+G^{-1}\star C^g_{+}-G^{-1}G^{-1}_{-}\Delta G\bigcdot\Delta C^r =G^{-1}\bigcdot C^r+G^{-1}\star C^g_{+}
\end{align*}
and
\begin{align*}
H&:=G^{-1}G^{-1}_{-}\Delta\xN\bigcdot A^o+G^{-1}_{-}\big(\wtG^{-1}-G^{-1}\big)\bigcdot [\xN,m]
\\&\quad -RG^{-1}\bigcdot A^o+Y_{-}G^{-1}G_{-}\wtG^{-1}\bigcdot A^o=\sum_{i=1}^4 H^i.
\end{align*}
We first recall that $\Delta A^o=\wtG-G$ and $\Delta m=\wtG-G_{-}$. Therefore,
\begin{align*}
&H^1+H^2=G^{-1}G^{-1}_{-}\Delta\xN\bigcdot A^o+G^{-1}_{-}\big(\wtG^{-1}-G^{-1}\big)\bigcdot [\xN,m] \\
&=G^{-1}G^{-1}_{-}\big(\Delta\xN\bigcdot A^o-\wtG^{-1}\Delta A^o \bigcdot [\xN,m] \big)=G^{-1}G^{-1}_{-}\big(\Delta\xN\bigcdot A^o-\wtG^{-1}\Delta \xN \Delta m\bigcdot A^o \big) \\
&=G^{-1}G^{-1}_{-}\big(\Delta\xN-\wtG^{-1}(\wtG-G_{-})\Delta\xN\big)\bigcdot A^o=G^{-1} \wtG^{-1}\Delta\xN\bigcdot A^o.
\end{align*}

Next, we deduce from \eqref{eq11.3} that
\[
\Delta \xN=\Delta (YG)+R\Delta A^o-\Delta C^r=\Delta (YG)+ R(\wtG-G )-\Delta C^r
\]
and thus
\begin{align*}
H&=G^{-1}\wtG^{-1}\big(\Delta (YG)+ R(\wtG-G)-\Delta C^r\big)\bigcdot A^o-R G^{-1}\bigcdot A^o +Y_{-}G^{-1}G_{-}\wtG^{-1}\bigcdot A^o\\
&=G^{-1} \wtG^{-1}\big(YG+ R(\wtG-G )-\Delta C^r\big)\bigcdot A^o -R G^{-1}\bigcdot A^o \\
&=\wtG^{-1}Y\bigcdot A^o-G^{-1}\wtG^{-1}\Delta C^r\bigcdot A^o-R\wtG^{-1}\bigcdot A^o \\
&=\wtG^{-1}Y\bigcdot A^o-G^{-1}\wtG^{-1}\Delta A^o\bigcdot C^r-R\wtG^{-1}\bigcdot A^o \\
&=(Y-\xR)\bigcdot \Gamma +(\wtG^{-1}-G^{-1})\bigcdot C^r.
\end{align*}
To complete the derivation of 
\eqref{eq11.6}, it suffices to substitute $K$ and $H$ into \eqref{eq11.8}.
\end{proof}

\end{document}